\documentclass[a4paper]{article}

\usepackage{amsfonts}
\usepackage{graphicx}
\usepackage{amsthm}
\usepackage{amsmath}
\usepackage{amssymb}
\usepackage{enumerate}
\allowdisplaybreaks

\usepackage{sagetex}

\newtheorem{proposition}{Proposition}[section]
\newtheorem{lemma}[proposition]{Lemma}
\newtheorem{corollary}[proposition]{Corollary}
\newtheorem{theorem}[proposition]{Theorem}
\newtheorem{definition}[proposition]{Definition}
\theoremstyle{definition}

\newtheorem{remark}[proposition]{Remark}
\numberwithin{equation}{section}

\newenvironment{mysage}{\sagesilent}{\endsagesilent}

\begin{document}

\begin{center}
\LARGE
\textbf{On short expressions for cosets of permutation subgroups}
\bigskip\bigskip

\large
Daniele Dona\footnote{The author was partially supported by the European Research Council under Programme H2020-EU.1.1., ERC Grant ID: 648329 (codename GRANT).}
\bigskip

\normalsize
Mathematisches Institut, Georg-August-Universit\"at G\"ottingen

Bunsenstra\ss e 3-5, 37073 G\"ottingen, Germany

\texttt{daniele.dona@mathematik.uni-goettingen.de}
\bigskip\bigskip\bigskip
\end{center}

\begin{minipage}{110mm}
\small
\textbf{Abstract.} Following Babai's algorithm \cite{Bab16a} for the string isomorphism problem, we determine that it is possible to write expressions of short length describing certain permutation cosets, including all permutation subgroups; this is feasible both in the original version of the algorithm and in its CFSG-free version, partially done by Babai \cite[\S 13.1]{Bab16a} and completed by Pyber \cite{Pyb16}. The existence of such descriptions gives a weak form of the Cameron-Mar\'oti classification even without assuming CFSG. We also thoroughly explicate Babai's recursion process (as given in Helfgott \cite{HBD17}) and obtain explicit constants for the runtime of the algorithm, both with and without the use of CFSG.
\medskip

\medskip

\end{minipage}
\bigskip

\begin{mysage}
### Rounding ####
def roundup(x,d):
    return round(ceil(x*10^d)/10^d,ndigits=d)
def rounddown(x,d):
    return round(floor(x*10^d)/10^d,ndigits=d)
### Redefining usual functions ###
# This will ensure that functions like exp interact well
# with the use of RIF(). Thanks are due to Marc Mezzarobba,
# who told it to Harald Helfgott who told it to me.
def expinterval(x):
    return x.exp()
def loginterval(x):
    return x.log()
def sqrtinterval(x):
    return expinterval(1/2*loginterval(x))
\end{mysage}

\section{Introduction}

Studying permutation subgroups is a rich part of today's research in finite group theory. The interest in permutations is even more understandable in light of the Classification of Finite Simple Groups (CFSG): a \textit{simple group} is a group that has no nontrivial normal subgroup, and simple groups are usually seen as the equivalent of prime numbers in group theory because of the Jordan-H\"{o}lder theorem \cite[\S 55]{Jor70} \cite{Hol89}. CFSG states that every finite simple group is either a cyclic group of size $p$ prime, $\mathrm{Alt}(n)$ for $n\geq 5$, a group of Lie type or one of $26$ exceptional groups (see for instance \cite[\S 1.2]{Wil09}).

CFSG has also many consequences, some of which we employ in the course of our reasoning. For example, it is possible to give better classification theorems of permutation subgroups using CFSG than not using it; on the other hand, while CFSG is generally accepted it is also very unwieldy, so that proving results without using CFSG is preferable to the alternative. Let us consider what we will call Theorem~\ref{th:cammar}, a consequence of a result by Cameron \cite{Cam81} and Mar\'oti \cite{Mar02} that describes all the primitive permutation groups as either having relatively small size or being very close to a wreath product of alternating groups: compare it with Pyber's result \cite{Pyb93} (Theorem~\ref{th:pyber}), that manages to give a similar description only for doubly transitive subgroups.

The theorem by Cameron and Mar\'oti is used in a recent result by Babai \cite{Bab16a} on the quasipolynomial procedure to solve the string and graph isomorphism problems (with \textit{quasipolynomial} we mean that it takes time $n^{O(\log^{c}n)}$, where $n$ is either the length of the strings or the size of the graphs involved and $c$ is some absolute constant): Cameron-Mar\'oti is the key passage to start the whole process and to keep the recursion running, and as we just said it depends on CFSG. However, it is possible to slightly modify Babai's proof to make it independent from CFSG: this modification process was initiated by Babai himself \cite[\S 13.1]{Bab16a} using Pyber's result; it was then completed by Pyber \cite{Pyb16} who proved what is called Lemma~4.1 in \cite{HBD17} without resorting to the Schreier conjecture, thus making Babai's algorithm CFSG-free at the price of making the bound worse (although still quasipolynomial, as Theorem~\ref{th:shortexpr} will show).

Our analysis here, on a first superficial level, provides a more explicit runtime for Babai's algorithm, both in the CFSG and the CFSG-free case. We will follow Helfgott's description of Babai's result given in \cite{Hel19} \cite{HBD17}, instead of Babai's original formulation in \cite{Bab16a}: Helfgott makes the algorithm more explicit and proves that the procedure actually takes time $n^{O(\log^{2}n)}$ when CFSG is available; we will make it even more explicit and determine the constants in front of the logarithm. Also, in \cite{Hel19} the reader's attention is justifiably focused on the proof of the single steps that are involved in the procedure, while the interstitial reasoning that details the recursion is only sketched: in \cite{Hel19}, this part is contained mostly in \S 3, \S 5.3, \S 6.2 and Appendix A; conversely, we will concentrate on the jumping between the main processes to delineate what the flow of the algorithm is, while using its individual theorems and subroutines as black boxes whose validity and well-functioning is taken for granted (we will mention the most important ones in \S\ref{se:majorou}). This will give us the control we need to determine the runtime with the desired accuracy.

On a deeper level, the way we achieve the goal described above is interesting on its own. Babai's algorithm is combinatorial in nature, although it is based on group-theoretic results; on the other hand, the combinatorial techniques developed by Babai have also been used before to deduce consequences for permutation subgroups, such as in \cite{Bab81}. It turns out that this is possible also in the case of Babai's quasipolynomial algorithm: since the procedure described by him is closely translatable to the CFSG-free case, it is possible to give a description of permutation subgroups that shares some characteristics of Cameron's result even when CFSG is not available, simply by making a subgroup pass through the algorithm, in a way that will be clarified in the next section; in brief, the use of the algorithm reveals structural information about permutation subgroups that we translate in the language of Theorem~\ref{th:shortexpr} as being able to write them as short expressions made of ``easy'' or ``atomic'' subgroups, where shortness here is just another face of the quasipolynomiality of the whole process.

That all of this can be useful, and that Theorem~\ref{th:shortexpr} can potentially do a job qualitatively similar to Cameron's theorem despite its different language, can be witnessed in \cite[\S 6]{Don20}. A decomposition similar to what we achieve in Theorem~\ref{th:shortexpr}, but based directly on Cameron, makes its appearance in \cite[Prop.~4.6]{Hel18} and is fundamental in proving a diameter bound for $\mathrm{Alt}(n)$ that goes through a sort of product theorem (i.e.\ a result like the key proposition in \cite{Hel08}). Passing through our decomposition instead, one would achieve the more modest and conditional result laid out in \cite[Thm.\ 6.3.6]{Don20}, which however shows already the potential power of our analysis.

%
%
%

\section{Standard definitions}\label{se:shortexprdef}

Before we start, let us recall here some standard terms and properties, coming from permutation group theory.

Let $A$ be a finite set: the \textit{symmetric group} $\mathrm{Sym}(A)$ is the group of all permutations of $A$, and any subgroup $G\leq\mathrm{Sym}(A)$ is called a \textit{permutation subgroup}; a particular permutation subgroup is the \textit{alternating group} $\mathrm{Alt}(A)$, which is the index $2$ subgroup that collects the even permutations of $A$ (i.e.\ the permutations obtainable as products of an even number of two-element transpositions). If $[n]=\{1,2,\ldots,n\}$, we write $\mathrm{Sym}(n),\mathrm{Alt}(n)$ for $\mathrm{Sym}([n]),\mathrm{Alt}([n])$.

\begin{definition}\label{de:transitive}
Let $n\geq 1$, and let $G\leq\mathrm{Sym}(n)$ be a permutation subgroup. $G$ is said to be transitive if for any two elements $x,y\in [n]$ there exists a $g\in G$ with $g(x)=y$. $G$ is intransitive if it is not transitive.

Let $d\geq 1$. $G$ is said to be $d$-transitive if for any two $d$-tuples of distinct elements $(x_{1},\ldots,x_{d}),(y_{1},\ldots,y_{d})\in [n]^{d}$ there is a $g\in G$ with $g(x_{i})=y_{i}$ for each $1\leq i\leq d$. A $2$-transitive subgroup is also referred to as doubly transitive.

The group $G\leq\mathrm{Sym}(n)$ is a giant if either $G=\mathrm{Sym}(n)$ or $G=\mathrm{Alt}(n)$.
\end{definition}

Transitive subgroups of $\mathrm{Sym}(n)$ have only one orbit for their natural action on $[n]$. There is another action of permutation subgroups that we will have to consider, namely the one on the set of $k$-subsets of $[n]$, denoted by $\binom{[n]}{k}$ (in obvious analogy with the binomial coefficients); in particular, the action of a $d$-transitive group on $\binom{[n]}{d}$ has only one orbit too. The same abstract group $G$ can be embedded into symmetric groups of different degrees, and thus be transitive or intransitive depending on the situation, therefore we will always specify ``$G\leq\mathrm{Sym}(n)$'' or similar notations to indicate that $G$ is considered to be of degree $n$; one of the reductions we operate, the one we call ``fourth action'' in \S\ref{se:shortexprcost}, is a passage to a smaller degree without changing $G$, so it is an important detail to keep in mind.

Let us see another important characteristic of the action of permutation groups.

\begin{definition}\label{de:primitive}
Let $G\leq\mathrm{Sym}(n)$ be transitive. A system of blocks of (the action of) $G$ is a partition $\mathcal{B}=\{B_{1},B_{2},\ldots,B_{r}\}$ of $[n]$ such that for every $g\in G$ and every $1\leq 1,j\leq r$ either $B_{i}=g(B_{j})$ or $B_{i}\cap g(B_{j})=\emptyset$. A trivial system of blocks is either the system $\mathcal{B}=\{[n]\}$ or the system $\mathcal{B}=\{\{1\},\{2\},\ldots,\{n\}\}$.

$G$ is primitive if the only systems of blocks it has are the trivial ones; $G$ is imprimitive if it is not primitive. $G$ is uniprimitive if it is primitive and not $2$-transitive.
\end{definition}

By transitivity, all the blocks of the same system have the same size. Every $2$-transitive group is primitive, but not vice versa: in other words, there exist uniprimitive groups, for example $\mathrm{Alt}(n)$ acting on $\binom{[n]}{2}$, provided that $n$ is large enough ($n=6$ is sufficient\footnote{Double transitivity fails because if $a$ is stabilized then all pairs containing $a$ are sent to each other; let us sketch the argument for primitivity. There is $1$ such that $\{1,2\}$ is in a block $B$ and $\{1,3\}$ is not: if there is $\{4,5\}\in B$, use $(2\ 3\ 6)$. If on the contrary all pairs in $B$ touch either $1$ or $2$, we can have either $\{1,3\}\in B$ and $\{1,4\}\not\in B$ (and use $(3\ 4\ 5)$) or every $\{1,x\},\{2,y\}\in B$ (and use $(1\ 2\ 3)$); if there are none at all, $|B|=1$.}). Similarly, there are transitive but imprimitive groups: an example of minimal size in terms of $|G|+n$ is $\langle(1\ 2),(1\ 3)(2\ 4)\rangle$ acting on $\{1,2,3,4\}$.

Finally, let us not miss an opportunity to describe the following action, since it plays a central role in Cameron.

\begin{definition}\label{de:wreath}
Let $G,H$ be finite groups acting on finite sets $V,W$ respectively. Then $G\wr H$, the wreath product of $G$ by $H$, is defined to be the semidirect product $G^{|W|}\rtimes H$; in other words, $G\wr H$ is the group whose underlying set is $G^{|W|}\times H$ and whose group operation is
\begin{equation*}
(g_{w_{1}},\ldots,g_{w_{|W|}},h)\cdot(g'_{w_{1}},\ldots,g'_{w_{|W|}},h')=(g_{w_{1}}g'_{h^{-1}(w_{1})},\ldots,g_{w_{|W|}}g'_{h^{-1}(w_{|W|})},hh').
\end{equation*}
The primitive action of $G\wr H$ on $V^{|W|}$ is defined to be
\begin{equation*}
(g_{w_{1}},\ldots,g_{w_{|W|}},h)\cdot(v_{w_{1}},\ldots,v_{w_{|W|}})=(g_{h^{-1}(w_{1})}v_{h^{-1}(w_{1})},\ldots,g_{h^{-1}(w_{|W|})}v_{h^{-1}(w_{|W|})}).
\end{equation*}
\end{definition}

There are several wreath products in more general contexts, but for us this will be sufficient. The primitive action of $G\wr H$ is also called product action or exponentiation in the literature \cite[\S 4.3]{Cam99} \cite[\S 2.7]{DM96} \cite[\S 4.1]{JK81}; there is also another natural action of the wreath product, the \textit{imprimitive action} on $V\times W$, but we will not encounter it.

In \cite[Thm.\ 1.1]{Mar02}, the definition above is used with $G=\mathrm{Sym}(m)$ and $H=\mathrm{Sym}(r)$ and their natural actions on $V=\binom{[m]}{k}$ and $W=[r]$ respectively.

\section{Main theorem: statement}\label{se:shortexprres}

Let us start with a permutation subgroup $G\leq\mathrm{Sym}(n)$. How ``easy'' is it to describe? Or rather, what are the ``easy'' permutation subgroups and how can we obtain all subgroups by building them out of the easy ones?

The easiest kind of subgroup that one can imagine would likely be a product of symmetric groups: given a partition $\{[n_{i}]\}_{i}$ of $[n]$, in the sense that $\sum_{i}n_{i}=n$, the subgroup corresponding to $\prod_{i}\mathrm{Sym}(n_{i})$ (provided that we fix a way to partition $[n]$ into these $[n_{i}]$) is very easily describable, in terms of generators, size, membership, etc...; we are curious about the way in which we can assemble groups of this sort to create $G$, or more generally a coset of $G$ if possible. Specifically, given a certain $H=\prod_{j}\mathrm{Sym}(n_{j})$ with $\sum_{j}n_{j}=n$ and a general $G\leq\mathrm{Sym}(n)$, we are going to give a description of cosets of the form $G\cap H\sigma$ in terms of easy subgroups; note that this does not include all the possible permutation cosets: for example, $G'\eta$ with $G'$ transitive is of the form $G\cap H\sigma$ only if $H=\mathrm{Sym}(n)$, which implies that $\eta$ is the identity permutation. On the other hand, by the same reasoning we promptly see that any subgroup $G'$ falls into this class of cosets. The reason why we restrict to these cosets will lie in our use of Babai's result (see Definition~\ref{de:iso}).

Let us define now more rigorously what it means to build an expression for $G\cap H\sigma$ starting from easy building blocks. Our \textit{atomic} elements are:
\begin{itemize}
\item[($\mathcal{A}$)]\label{eq:atom} cosets $G\sigma$ of permutation subgroups $G$ of the form $\mathrm{Alt}(\bigsqcup_{i}A_{i})\cap\prod_{i}\mathrm{Sym}(A_{i})$ (where the $A_{i}$ are disjoint sets).
\end{itemize}
So the atoms are defined to be the cosets of the even permutation part of the aforementioned ``easiest subgroups''. In particular, the trivial subgroup $\{\mathrm{Id}_{|\Omega|}\}$ is an atom, being simply $\mathrm{Sym}(1)^{|\Omega|}$, and so are all singletons $\{\sigma\}$, being its cosets.

We declare the atoms to be \textit{well-formed}. We can combine well-formed expressions to form more complex ones; the legitimate ways to do it are the following three.
\begin{enumerate}[($\mathcal{C}$1)]
\item\label{eq:combunion} \textit{Paste cosets of a subgroup to get the whole group.}

Let $G'\leq G\leq\mathrm{Sym}(A)$ with $\{\sigma_{i}\}_{i}$ a set of representatives of $G'$ in $G$, and let $H=\prod_{j}\mathrm{Sym}(A_{j})$ for some partition $\{A_{j}\}_{j}$ of $A$; suppose that for some fixed $\sigma\in\mathrm{Sym}(A)$ the cosets $G'\cap H\sigma\sigma_{i}^{-1}$ are all well-formed: then $G\cap H\sigma=\bigcup_{i}(G'\cap H\sigma\sigma_{i}^{-1})\sigma_{i}$ is also well-formed.

\item\label{eq:combglue} \textit{Paste disjoint domains to get a group acting on both.}

Let $G\leq\mathrm{Sym}(A_{1})\times\mathrm{Sym}(A_{2})$; for $i=1,2$, let $\pi_{i}:G\rightarrow\mathrm{Sym}(A_{i})$ be the natural projections, let $H_{i}=\prod_{j}\mathrm{Sym}(A_{ij})$ for some partition $\{A_{ij}\}_{j}$ of $A_{i}$, and let $\sigma_{i}\in\mathrm{Sym}(A_{i})$. Suppose that $\pi_{1}(G)\cap H_{1}\sigma_{1}=K\tau$ is well-formed, and suppose that $\pi_{2}(\pi_{1}^{-1}(K))\cap H_{2}\sigma_{2}\pi_{2}(\pi_{1}^{-1}(\tau))$ is well-formed too: then $G\cap (H_{1}\times H_{2})(\sigma_{1},\sigma_{2})$ is well-formed.

\item\label{eq:combalt} \textit{Paste a group fixing a set of blocks with an alternating group permuting them.}

Let $G\leq\mathrm{Sym}(A)$ be a well-formed subgroup, contained in $\prod_{i}\mathrm{Sym}(A_{i})$ for some partition $\{A_{i}\}_{i}$ of $A$ into equally sized parts; let $\sigma_{1},\sigma_{2},\sigma'$ be three permutations of $A$ and suppose that $\langle\{\sigma_{1},\sigma_{2}\}\rangle$ permutes the $A_{i}$ in the same way as $\mathrm{Alt}(\Gamma)$ permutes $\binom{\Gamma}{k}$ for some $\Gamma,k$: then $\langle G\cup\{\sigma_{1},\sigma_{2}\}\rangle\sigma'$ is also well-formed.
\end{enumerate}

Since the trivial subgroup is an atom, all subgroups $G$ could be written as a well-formed expression by ($\mathcal{C}$1), choosing $G'=\{\mathrm{Id}_{|\Omega|}\}$, $H=\mathrm{Sym}(\Omega)$ and any $\sigma$. That is uninteresting, though, since we need $|G|$ atoms to perform such a task: the point is to use as few of them as possible. Our main theorem gives a way to build a well-formed expression of small length for $G$, and even for any $G\cap H\sigma$.

\begin{mysage}
main_cfsg=103
main_free=26
\end{mysage}

\begin{theorem}\label{th:shortexpr}
Let $n\geq 1$, let $G\leq\mathrm{Sym}(n)$, let $H=\prod_{i}\mathrm{Sym}(\Sigma_{i})$ for some partition $\{\Sigma_{i}\}_{i}$ of $[n]$, and let $\sigma\in\mathrm{Sym}(n)$.

Then, we can write a well-formed expression for $G\cap H\sigma$, starting from atomic elements $\mathrm{(\mathcal{A})}$ and combining them using $\mathrm{(\mathcal{C}1)}$-$\mathrm{(\mathcal{C}2)}$-$\mathrm{(\mathcal{C}3)}$, such that the number of atomic elements involved in the construction is bounded by $n^{K\log^{c}n}$, where $(K,c)=(\sage{main_cfsg},2)$ if we assume CFSG and $(K,c)=(\sage{main_free}e^{1/\varepsilon^{2}},5+\varepsilon)$ otherwise for any $\varepsilon>0$ small enough.

The time necessary to find such an expression is bounded by $O(n^{11+K\log^{c}n})$.
\end{theorem}

One can verify that $\varepsilon<\frac{1}{100}$ is indeed small enough.

The similarities with \cite[Prop.~4.6]{Hel18} are important, as they are exactly of the nature that we would need to free the bound on $\mathrm{diam}(\mathrm{Alt}(n))$ proved therein from the use of CFSG: the descent to smaller cosets (or ascent to larger ones, for us) works in the same way, and the quasipolynomial bound is fundamental for the diameter. The only difference that prevents a direct substitution is the fact that $\mathrm{(\mathcal{C}1)}$ allows for any subgroup, instead of restricting to normal subgroups like we would need for other procedures given in the course of such a proof. See \cite[\S 6]{Don20} for a more in-depth analysis of this point.

The runtime claimed in Theorem~\ref{th:shortexpr} is in reality a bound on the runtime for Babai's algorithm: the construction process of the well-formed expression, as illustrated in the following sections, is part of the description process necessary to solve the string isomorphism problem; in the proof we will calculate the cost for the latter, thus retrieving a bound for the former as well.

Setting aside the time issue, this theorem does not surprise us if we assume CFSG. Cameron implies in its stronger form that any primitive permutation subgroup either is small enough to be expressed as the union of $\leq n^{O(\log^{2}n)}$ singletons through $\mathrm{(\mathcal{C}1)}$ or it has as large subgroup a wreath product $\mathrm{Alt}(\Gamma)\wr\mathrm{Alt}(s)$ where $\mathrm{Alt}(\Gamma)$ acts on $\binom{\Gamma}{k}$ (see Definition~\ref{de:wreath}), so that it is susceptible of being described using repeatedly $\mathrm{(\mathcal{C}3)}$; if the subgroup is not primitive, it is not difficult to reduce to this case by working on each block separately and then uniting and glueing together the pieces with $\mathrm{(\mathcal{C}1)}$ and $\mathrm{(\mathcal{C}2)}$.

Without assuming CFSG however, the situation changes. It is true that, for doubly transitive permutation subgroups, Theorem~\ref{th:shortexpr} would be a consequence of Pyber's result: either such a group is $\mathrm{Sym}(n)$ or $\mathrm{Alt}(n)$, or it has size $\leq n^{O(\log^{2}n)}$; the discussion goes basically as above. Pyber's result does not however say anything about subgroups that are transitive but not doubly transitive; in this sense, Theorem~\ref{th:shortexpr} extends this CFSG-free description to this class of permutation subgroups as well (and \cite[Prop.~4.6]{Hel18} is needed for all transitive groups).

One last note: the computation of $K$ in the main theorem, and many of the intermediate results leading to it, have been performed with SageMath, version 8.9. The calculations are elementary enough to be easily reproducible with any software, but SageMath is open-source and can be embedded into LaTeX, which is why the author chose to use it.

\section{Elementary routines}\label{se:shortexprback}

Let us define the fundamental objects in the study of SIP.

\begin{definition}\label{de:iso}
Let $\Omega$ be a finite set, let $G\leq\mathrm{Sym}(\Omega)$ and let $\mathbf{x},\mathbf{y}:\Omega\rightarrow\Sigma$ be two strings. The set of isomorphisms from $\mathbf{x}$ to $\mathbf{y}$ in $G$ is defined as
\begin{equation*}
\mathrm{Iso}_{G}(\mathbf{x},\mathbf{y})=\{g\in G|\mathbf{x}^{g}=\mathbf{y}\}=\{g\in G|\forall r\in\Omega(\mathbf{x}(r)=\mathbf{y}(g(r)))\}.
\end{equation*}
The group of automorphisms of $\mathbf{x}$ in $G$ is defined as $\mathrm{Aut}_{G}(\mathbf{x})=\mathrm{Iso}_{G}(\mathbf{x},\mathbf{x})$.
\end{definition}

The sets of isomorphisms $\mathrm{Iso}_{G}(\mathbf{x},\mathbf{y})$ are precisely the intersections $G\cap H\sigma$, $H$ being a product of smaller symmetric groups, that are featured in Theorem~\ref{th:shortexpr}: in fact, a permutation of $\Omega$ is in such a set if and only if it is in $G$ and for every letter of $\Sigma$ it sends the preimage of that letter in $\mathbf{x}$ to its preimage in $\mathbf{y}$. $H$ is therefore $\prod_{\alpha\in\mathbf{x}(\Omega)}\mathrm{Sym}(\mathbf{x}^{-1}(\alpha))$, and vice versa, given a product of symmetric groups and a $\sigma$, it is possible to define $\mathbf{x}$ as being piecewise constant with a letter for each symmetric group and then define $\mathbf{y}=\mathbf{x}^{\sigma}$.

This also reveals how to find an expression for any permutation subgroup $G\leq\mathrm{Sym}(\Omega)$: this corresponds to finding $\mathrm{Aut}_{G}(\alpha^{|\Omega|})$, where $\alpha^{|\Omega|}$ is the constant string consisting of one letter repeated $|\Omega|$ times, or in other words to making the algorithm run ``in neutral'' on a trivial string so as to capture only $G$.

\begin{remark}\label{re:iso}
Every time we describe $\mathrm{Iso}_{G}(\mathbf{x},\mathbf{y})$ as a coset $G'\tau$, where $G'\leq\mathrm{Sym}(\Omega)$ and $\tau\in\mathrm{Sym}(\Omega)$, $G'$ is actually $\mathrm{Aut}_{G}(\mathbf{x})$ and $\tau$ is an element of $G$ sending $\mathbf{x}$ to $\mathbf{y}$.

In fact, since $G'$ is a subgroup of $\mathrm{Sym}(\Omega)$ it contains the trivial permutation, so that $\tau\in\mathrm{Iso}_{G}(\mathbf{x},\mathbf{y})$: this proves what we claimed about $\tau$. If $g\in G'$ (so that $g\tau$ sends $\mathbf{x}$ to $\mathbf{y}$) then $g$ fixes $\mathbf{x}$ since permutations are bijections and any $\mathbf{x}'\neq\mathbf{x}$ will not be sent to $\mathbf{y}$ by $\tau$; therefore by definition $g$ is also an element of $\mathrm{Aut}_{G}(\mathbf{x})$. On the other hand, if $\sigma\in\mathrm{Aut}_{G}(\mathbf{x})$ then $\sigma\tau\in\mathrm{Iso}_{G}(\mathbf{x},\mathbf{y})=G'\tau$ and $\sigma\in G'$; this proves also that $G'=\mathrm{Aut}_{G}(\mathbf{x})$.
\end{remark}

We begin by providing several simple results on computations that we have to constantly perform throughout the whole procedure. Before that, a couple of definitions; if $G\leq\mathrm{Sym}(\Omega)$ and $\Delta\subset\Omega$, the \textit{setwise stabilizer} and the \textit{pointwise stabilizer} of $\Delta$ are respectively
\begin{align*}
G_{\Delta} & = \{g\in G|g(\Delta)=\Delta\}, \\
G_{(\Delta)} & = \{g\in G|\forall r\in\Delta(g(r)=r)\}.
\end{align*}
We also write $G_{(r_{1},\ldots,r_{i})}$ for $G_{(\{r_{1},\ldots,r_{i}\})}$. Trying to find the setwise stabilizer for a generic $\Delta$ is a task of difficulty comparable to producing $\mathrm{Iso}_{G}(\mathbf{x},\mathbf{y})$ itself; on the other hand, producing pointwise stabilizers is much easier (see Corollary~\ref{co:schreier}\eqref{co:schreierpwstab}), and we can walk down this route to obtain basic but useful algorithms.

\begin{proposition}[Schreier-Sims algorithm]\label{pr:schreier}
Let $\Omega=\{x_{1},x_{2},\ldots,x_{n}\}$ and let $G\leq\mathrm{Sym}(\Omega)$ be provided with a set of generators $A$. Then there is an algorithm that finds in time $O(n^{5}+n^{3}|A|)$ a set $C$ of generators of $G$ of size $\leq n^{2}$ such that for every $0\leq i\leq n-2$ and for every coset of $G_{(x_{1},\ldots,x_{i},x_{i+1})}$ inside $G_{(x_{1},\ldots,x_{i})}$ there exists a unique $\gamma\in C$ that is a representative of that coset.
\end{proposition}

\begin{proof}
See \cite[\S 1.2]{Luk82} or \cite[Alg.~1]{Hel19}.
\end{proof}

We will see that in our base cases corresponding to the atoms ($\mathcal{A}$) the number of generators will be polynomial in $n$, so that we will not have problems supposing that the Schreier-Sims algorithm takes polynomial time in $n$; from now on, when we talk about polynomial time (or size, or cost) we mean polynomial in $n$, the length of the strings involved. It also happens at some point that we take the union of several cosets, and the process produces sets of generators of size comparable to the number of cosets (as described in Proposition~\ref{pr:Gimprim}); in that case, the time will be more conspicuous: for instance, Corollary~\ref{co:cammar}\eqref{co:cammarsmall} and Proposition~\ref{pr:pybersmall} entail a cost of order $m^{O(\log^{2}n)}n^{O(1)}$ for the filtering of generators through Schreier-Sims.

In any case, every time a $G$ is already ``given'', or has been ``described'' or ``determined'', or other similar locutions, we will suppose that it has a quadratic number of generators thanks to Schreier-Sims (unless explicitly stated otherwise).

Proposition~\ref{pr:schreier} provides us with many useful polynomial-time procedures, as shown below.

\begin{corollary}\label{co:schreier}
Let $|\Omega|=n$ and let $G\leq\mathrm{Sym}(\Omega)$ be provided with a set of generators $A$ of polynomial size. Then the following tasks can be accomplished in polynomial time:
\begin{enumerate}[(a)]
\item\label{co:schreiersize} determine $|G|$;
\item\label{co:schreierbelongs} determine whether a certain $g\in\mathrm{Sym}(\Omega)$ is in $G$;
\item\label{co:schreiersubtest} given a subgroup $H\leq G$ with index $[G:H]$ of polynomial size and given a polynomial-time test that determines whether a certain $g\in G$ is in $H$, determine $H$ and a representative of each coset of $H$ in $G$;
\item\label{co:schreierpreim} given a homomorphism $\varphi: G\rightarrow\mathrm{Sym}(\Omega')$ with $\Omega'$ of polynomial size and given a subgroup $H\leq\mathrm{Sym}(\Omega')$, determine $\varphi^{-1}(H)$, or given an element $\tau\in\mathrm{Sym}(\Omega')$, determine an element of $\varphi^{-1}(\tau)$;
\item\label{co:schreierpwstab} given a set $S\subseteq\Omega$, determine $G_{(S)}$;
\item\label{co:schreiersetstab} provided that $G$ acts transitively imprimitively on $\Omega$ and given a system of blocks of its action on $\Omega$, determine the stabilizer of this system;
\end{enumerate}
Moreover, we can explicitly write in time $O(n^{5}+n^{3}|A|+n^{2}|G|)$ all the elements of $G$.
\end{corollary}

\begin{proof}
For parts \eqref{co:schreiersize}-\eqref{co:schreierbelongs}-\eqref{co:schreiersubtest} see \cite[Ex.~2.1(a)-2.1(c)]{Hel19}, based on \cite[Cor.~1]{FHL80} and \cite[Lemma~1.2]{Luk82}; the representatives in part \eqref{co:schreiersubtest} are the elements of $C_{-1}$ in the solution of \cite[Ex.~2.1(c)]{Hel19} given in \cite[App.~B]{HBD17}\footnote{Between \cite{Hel19} and \cite{HBD17}, Exercise 2.1(b) in one corresponds to Exercise 2.1(c) in the other. The author apologizes, but that was the order in which he proved things during the translation process: if he had respected the original order, part (b) would have depended on part (c).}. Part \eqref{co:schreierpreim} is similar to \eqref{co:schreiersubtest}, see \cite[Ex.~2.1(b)]{Hel19}; finding an element of the preimage of a generator is a passage inside the proof of the procedure that finds $\varphi^{-1}(H)$, so to solve the second issue we can take $H=\langle\tau\rangle$. Finding pointwise stabilizers $G_{(S)}$ is a byproduct of Schreier-Sims itself, so we simply have to order $\Omega$ so that $S=\{x_{1},\ldots,x_{|S|}\}$ and Proposition~\ref{pr:schreier} will solve part \eqref{co:schreierpwstab} directly. Part \eqref{co:schreiersetstab} is an application of \eqref{co:schreierpreim}: $\Omega'$ will be the system of blocks (which means that $|\Omega'|<n$) and $H=\{\mathrm{Id}_{|\Omega'|}\}$.

The last statement is a consequence of the particular structure of the set of generators $C$ found through Schreier-Sims: $C$ is divided into sets $C_{0},\ldots,C_{n-2}$, each consisting of the generators $\gamma\in G_{(x_{1},\ldots,x_{i})}\setminus G_{(x_{1},\ldots,x_{i+1})}$, and each element of $G$ is written uniquely as a product $\gamma_{0}\gamma_{1}\ldots\gamma_{n-2}$ with $\gamma_{i}\in C_{i}$. There are $|G|$ such products, and a product of ttwo permutations is computable in time $O(n)$, whence the result.

Let us include here the runtimes of the other items, too. Parts \eqref{co:schreiersize}-\eqref{co:schreierbelongs}-\eqref{co:schreierpwstab} consist in using the Schreier-Sims algorithm at most twice with at most one more generator, so the runtime is $O(n^{5}+n^{3}|A|)$. In Schreier-Sims, the time is more explicitly of order $n\cdot(n^{2}\cdot n^{2}+n^{2}\cdot|A|)$, where $n$ comes from the use of the subroutine \textsc{Filter} in \cite[Alg.~1]{Hel19} and $n^{2}$ is the bound on the size of the final $C$; by this analysis, part \eqref{co:schreiersubtest} employs time $O(n^{2i+t}+n^{i+t}|A|)$, where $i$ is the maximum between $2$ and the exponent of the index $[G:H]$ and $t$ is the maximum between $1$ and the exponent of the test time for $H$. For part \eqref{co:schreierpreim}, we use Schreier-Sims first on $G$, then on each preimage of $\mathrm{Sym}(\Omega')_{(x'_{1},\ldots,x'_{i})}$, then we express each generator of $H$ as product of images of generators of $G$: this takes time $O(n^{5s}+n^{3}|A|+n^{h+2s})$, where $s$ is the maximum between $2$ and the exponent of $|\Omega'|$ and $h$ is the exponent of the number of generators of $H$. Using \eqref{co:schreierpreim}, part \eqref{co:schreiersetstab} takes time $O(n^{10}+n^{3}|A|)$.
\end{proof}

All these polynomial costs will not be particularly relevant: in the course of our reasoning we will not encounter an exponent of a polynomial cost that is larger than $14$, and this is negligible against the $n^{K\log^{e}n}$ we have at the end. The constants hidden in the big O notation are only depending on the cost of procedures like reading, writing, comparing elements, etc...: we will not care about them, but just carry them around inside the O.

Another important polynomial-time algorithm is the one illustrated in the following lemma: recalling the definition of transitivity and primitivity for permutation subgroups, it is clear that being able to quickly determine respectively orbits and blocks of the actions of groups that do not present these two properties is a beneficial skill for us to possess.

\begin{lemma}\label{le:orbblo}
Let $|\Omega|=n$ and $G\leq\mathrm{Sym}(\Omega)$. Then the orbits of the action of $G$ on $\Omega$ can be determined in time $O(n^{3})$; also, if $G$ is transitive but imprimitive, a system of minimal blocks for the action of $G$ on $\Omega$ can be determined in time $O(n^{4})$.
\end{lemma}

\begin{proof}
To determine the orbits, we follow \cite[Ex.~B.2]{HBD17}. Let $A$ be a set of generators of $G$, which by Schreier-Sims we can suppose is of size $\leq n^{2}$: the sets $A_{x}=\{x^{a}|a\in A\}$ for every $x\in\Omega$ can be determined in time $O(n^{3})$. After that, we follow this procedure: we start with any fixed $x_{0}\in\Omega$ and set $\Delta_{x_{0}}=\{x_{0}\}\cup A_{x_{0}}$; we divide the elements of $\Delta_{x_{0}}$ in ``examined'' (at this stage, only $x_{0}$) and ``unexamined'' (the other elements of $\Delta_{x_{0}}$). Then at every step we take an unexamined $x\in \Delta_{x_{0}}$ and we update $\Delta_{x_{0}}$ by adding the elements of $A_{x}$ to it: the newly added elements are marked as unexamined, while $x$ now is examined; the procedure stops when $\Delta_{x_{0}}$ becomes the orbit $\{x_{0}^g|g\in G\}$. If there is an element $x_{1}$ that has not yet been considered, we define $\Delta_{x_{1}}=\{x_{1}\}\cup A_{x_{1}}$ and go through the whole procedure again, until we have considered all the elements of $\Omega$: the final sets $\Delta_{x_{0}},\Delta_{x_{1}},\ldots,\Delta_{x_{m}}$ are the orbits of the action of $G$ on $\Omega$; this part takes time $O(n)$, so the runtime of the whole algorithm is $O(n^{3})$.

Suppose now that $G$ is transitive imprimitive: to determine the blocks we follow \cite[\S 2.1.2]{Hel19}, which is based on an idea by Higman (through Sims and then Luks). The idea in the previous case was basically to follow the edges of the Schreier graph of $G$ with set of generators $A$ on $\Omega$: we will do the same with different graphs now. Our preparatory work this time consists in considering all the pairs $\{x,x'\}\subseteq\Omega$ and constructing the sets $A_{x,x'}=\{\{x^{a},x'^{a}\}|a\in A\}$ in time $O(n^{4})$, forming a first graph; then we fix $x_{0}\in\Omega$ and for every other $x\in\Omega$ we build the following graph: the set of vertices is $\Omega$ and the edges are the pairs contained in the connected component of $\{x_{0},x\}$ of the first graph (finding the connected component takes linear time in the number of vertices, so $O(n^{2})$ here). In the newly formed graphs, the connected components containing $\{x_{0},x\}$ are the smallest blocks containing $\{x_{0},x\}$ (see \cite[Prop.~4.4]{Sim67}; again, finding the connected components is a $O(n)$ routine): once we find among the blocks constructed from each $x$ a block that is properly contained in $\Omega$, which exists for $G$ imprimitive, we can find a whole system by taking the other components of the graph given by the same $x$. The system may not be minimal, but we have only to repeat the whole process working with the set of blocks instead of $\Omega$; since at each iteration the blocks are at least twice the size of the ones at the previous step, eventually we reach a system that has blocks of maximal size, i.e.\ a minimal system. The whole process works in time $O\left(n^{4}+\left(\frac{n}{2}\right)^{4}+\left(\frac{n}{2^{2}}\right)^{4}+\ldots\right)=O(n^{4})$.
\end{proof}

Finally, we illustrate several equalities among different sets of isomorphisms (employed here in a slightly more flexible way than Definition~\ref{de:iso}) that will allow us to pass from difficult problems to easier ones, or to break down problems into smaller ones.

\begin{lemma}\label{le:chain}
Let $|\Omega|=n$, $G\leq\mathrm{Sym}(\Omega)$, $\sigma\in\mathrm{Sym}(\Omega)$ and let $\mathbf{x},\mathbf{y}:\Omega\rightarrow\Sigma$ be two strings. For $\Delta\subseteq\Omega$ invariant under $G,\sigma$, define the set of partial isomorphisms $\mathrm{Iso}^{\Delta}_{G\sigma}(\mathbf{x},\mathbf{y})$ as in Definition~\ref{de:iso} with $g\in G\sigma$ and $\mathbf{x}(r)=\mathbf{y}(g(r))$ necessary only for $r\in\Delta$.
\begin{enumerate}[(a)]
\item\label{le:chain1} We can pass from cosets to groups using
\begin{equation*}
\mathrm{Iso}^{\Delta}_{G\sigma}(\mathbf{x},\mathbf{y})=\mathrm{Iso}^{\Delta}_{G}(\mathbf{x},\mathbf{y}^{\sigma^{-1}})\sigma.
\end{equation*}
\item\label{le:chain2} We can split unions of cosets using
\begin{equation*}
\mathrm{Iso}^{\Delta}_{G\sigma_{1}\cup G\sigma_{2}}(\mathbf{x},\mathbf{y})=\mathrm{Iso}^{\Delta}_{G\sigma_{1}}(\mathbf{x},\mathbf{y})\cup\mathrm{Iso}^{\Delta}_{G\sigma_{2}}(\mathbf{x},\mathbf{y}).
\end{equation*}
\item\label{le:chain3} We can split unions of windows using
\begin{equation*}
\mathrm{Iso}^{\Delta_{1}\cup\Delta_{2}}_{G\sigma}(\mathbf{x},\mathbf{y})=\mathrm{Iso}^{\Delta_{2}}_{G_{1}}(\mathbf{x},\mathbf{y}^{\sigma_{1}^{-1}})\sigma_{1},
\end{equation*}
where $\mathrm{Iso}^{\Delta_{1}}_{G\sigma}(\mathbf{x},\mathbf{y})=G_{1}\sigma_{1}$.
\item\label{le:chain4} For every $g\in G$, call $g|_{\Delta}$ its restriction to $\Delta$, defined by simply forgetting what happens in $\Omega\setminus\Delta$ (since $G$ leaves $\Delta$ invariant, this is well-defined); define $S|_{\Delta},H|_{\Delta},\mathbf{x}|_{\Delta}$ for any $S\subseteq G$, $H\leq G$, $\mathbf{x}:\Omega\rightarrow\Sigma$ analogously. For any $h\in G|_{\Delta}$, let $\overline{h}$ be any element of $G$ whose restriction to $G|_{\Delta}$ is $h$; if $H\leq G|_{\Delta}$, define $\overline{H}$ analogously as the subgroup of $G$ whose restriction to $G|_{\Delta}$ is $H$ (since $G$ leaves $\Delta$ invariant, $\overline{H}$ is indeed a subgroup).

We can eliminate windows using
\begin{equation*}
\mathrm{Iso}_{G}^{\Delta}(\mathbf{x},\mathbf{y})=\overline{G'}\overline{\sigma},
\end{equation*}
where $\mathrm{Iso}_{G|_{\Delta}}(\mathbf{x}|_{\Delta},\mathbf{y}|_{\Delta})=G'\sigma$; this is independent from the choice of $\overline{\sigma}$.
\end{enumerate}
\end{lemma}

\begin{proof}
\eqref{le:chain1} It is easy from the definition: inside $\Delta$, the permutation $g=g'\sigma\in G\sigma$ sends $\mathbf{x}$ to $\mathbf{y}$ if and only if $g'$ sends $\mathbf{x}^{\sigma}$ to $\mathbf{y}$, i.e.\ if and only if it sends $\mathbf{x}$ to $\mathbf{y}^{\sigma^{-1}}$.

\eqref{le:chain2} It is obvious from the definition, since both sides mean the exact same thing, allowing in both cases $g$ to be either in $G\sigma_{1}$ or in $G\sigma_{2}$.

\eqref{le:chain3} First, we obtain $\mathrm{Iso}^{\Delta_{1}\cup\Delta_{2}}_{G\sigma}(\mathbf{x},\mathbf{y})=\mathrm{Iso}^{\Delta_{2}}_{G_{1}\sigma_{1}}(\mathbf{x},\mathbf{y})$ easily by examining the definitions: both sides simply mean that $g\in G\sigma$ has to respect both windows $\Delta_{1},\Delta_{2}$. Then we get $\mathrm{Iso}^{\Delta_{2}}_{G_{1}\sigma_{1}}(\mathbf{x},\mathbf{y})=\mathrm{Iso}^{\Delta_{2}}_{G_{1}}(\mathbf{x},\mathbf{y}^{\sigma_{1}^{-1}})\sigma_{1}$ from part \eqref{le:chain1}.

\eqref{le:chain4} $G'\sigma$ is the collection of permutations of $\Delta$ that send $\mathbf{x}$ to $\mathbf{y}$ as far as $\Delta$ is able to perceive. Passing to the whole $\Omega$ by considering $\overline{G'}$ and $\overline{\sigma}$, the result is the definition itself of $\mathrm{Iso}_{G}^{\Delta}(\mathbf{x},\mathbf{y})$.
\end{proof}

\begin{mysage}
Clim=102
mlim=8308
Clim_free=25
mlim_free=max(0,Clim_free)
\end{mysage}

\begin{remark}\label{re:thism}
In the future we are going to need to differentiate the cases of $n$ large and $n$ small. This will come in the form of $C\log^{c}n\leq n$, for certain $C,c>0$: if such an inequality is true, which would allow us to have an intermediate integer $m$ between them when needed, then $n$ is considered large. Let us make now this choice.

Assuming CFSG, we suppose that largeness means $\sage{Clim}\log^{2}n<m\leq n$, which implies $m,n\geq\sage{mlim}$. See \eqref{eq:multfinal} inside the proof of the main theorem, which is the final quantity to optimize. Without assuming CFSG we suppose instead that largeness means $\sage{Clim_free}e^{1/\varepsilon^{2}}(\log n)^{4+\varepsilon}<m\leq n$, which implies in particular $m,n\geq\sage{mlim_free}e^{1/\varepsilon^{2}}$. For $\varepsilon$ small (say $\varepsilon<\frac{1}{10}$), the CFSG-free condition is a stronger restriction.
\end{remark}

\section{Major routines}\label{se:majorou}

Before we turn to the algorithm itself, let us describe separately a couple of major routines that were introduced for the first time by Babai. We will not prove their validity here: both \cite{Bab16a} and \cite{Hel19} do that for us. What we want is to sum up their contribution to the runtime.

We start with a theoretical result, needed to differentiate between the CFSG and the CFSG-free case.

\begin{lemma}\label{le:epicfsg}
Let $G\leq\mathrm{Sym}(n)$ be primitive, and let $\phi:G\rightarrow\mathrm{Alt}(a)$ be an epimorphism.
\begin{enumerate}[(a)]
\item\label{le:epicfsgbabai} Assuming CFSG, if $a>\max\{8,2+\log_{2}n\}$ then $\phi$ is an isomorphism.
\item\label{le:epicfsgpyber} Not assuming CFSG, if $a>\max\{e^{1/\varepsilon^{2}},(\log_{2}n)^{4+\varepsilon}\}$ then $\phi$ is an isomorphism, for any $\varepsilon>0$ small enough.
\end{enumerate}
\end{lemma}

\begin{proof}
For \eqref{le:epicfsgbabai} see \cite[Lemma 8.2.4]{Bab16a} or \cite[Lemme 4.1]{Hel19}. For \eqref{le:epicfsgpyber} see \cite[Lemma 12]{Pyb16}, which states that it's sufficient to take $a>\max\{C,\log_{2}^{5}n\}$ for some constant $C$ (his ``$\log$'' is our ``$\log_{2}$''). Let us compute our version of the bound.

Using \cite[Thm.~C(v)]{MR96}, the sum of the first $s$ primes for $s\geq 6$ is bounded by $\frac{1}{2}s^{2}(\log s+\log_{2}s)$, so $\mathrm{Alt}(a)$ contains a cyclic subgroup whose order is the product of the first $\left\lceil\frac{5}{4}\sqrt{\frac{a}{\log a}}\right\rceil$ primes for $a\geq 10000$ (say). From \cite[Thm.~7]{Pyb16} and $a>(\log_{2}n)^{4+\varepsilon}$, if $\phi$ were not an isomorphism we would get $\left\lceil\frac{5}{4}\sqrt{\frac{a}{\log a}}\right\rceil<2a^{\frac{2}{4+\varepsilon}}<2a^{\frac{1}{2}-\frac{\varepsilon}{10}}$ (for small $\varepsilon$), which can be true only if $a\leq e^{1/\varepsilon^{2}}$ (again for small $\varepsilon$).
\end{proof}

A short verification shows that $\varepsilon<\frac{1}{100}$ is plenty enough for the result above to hold.

We are using Lemma~\ref{le:epicfsg} in the computation of the runtime of the following routine. The production and aggregation of \textit{local certificates} (see \cite[\S 10]{Bab16a} or \cite[\S 6]{Hel19}) is an important part of the algorithm.

\begin{mysage}
alim_low=roundup( 1/loginterval(RIF(2))+2/loginterval(RIF(mlim)) ,5)
alim=roundup( alim_low+1.00001/loginterval(RIF(mlim)) ,5) #Must be at most Clim/4.
alim_free=rounddown( Clim_free/4 ,5)
alim_freelow=rounddown( alim_free-1.00001/(expinterval(30)*loginterval(RIF(mlim_free))^4) ,5) #Must be larger than 1/(log 2)^5.
#alim_low=round(  ,ndigits=5)
#alim=round( 1.73036 ,ndigits=5)
#alim_freelow=round(  ,ndigits=5)
#alim_free=round( 1.00001 ,ndigits=5)
\end{mysage}

\begin{proposition}\label{pr:computcert}
Let $G\leq\mathrm{Sym}(n)$, and let $\phi:G\rightarrow\mathrm{Alt}(m)$ be an epimorphism; let $\mathbf{x}$ be a string of length $n$. Then we can find the group $F\leq\mathrm{Aut}_{G}(\mathbf{x})$ generated by the certificates of fullness in the time taken by $\frac{1}{2}m^{2a}naa!$ calls of the whole algorithm for strings of length $\leq\frac{n}{a}$, where
\begin{enumerate}[(a)]
\item\label{pr:computcertcfsg} $a\in\left(\sage{alim_low},\sage{alim}\right)\cdot\log n$, for $\sage{Clim}\log^{2}n<m\leq n$ (assuming CFSG), or
\item\label{pr:computcertfree} $a\in(\sage{alim_freelow},\sage{alim_free})\cdot e^{1/\varepsilon^{2}}(\log n)^{4+\varepsilon}$, for $\sage{Clim_free}e^{1/\varepsilon^{2}}(\log n)^{4+\varepsilon}<m\leq n$ for any $\varepsilon>0$ small enough (without assuming CFSG),
\end{enumerate}
and in both cases $a\leq\frac{m}{4}$, plus some additional time $O(m^{2a}n^{11})$.
\end{proposition}

\begin{proof}
The proof is contained in \cite[\S 6.1]{Hel19}. We are going to discuss the details of the runtime.

Let $T,T'$ be two ordered $a$-tuples of elements of $\Gamma$, where $a$ is as in Lemma~\ref{le:epicfsg}. Updating one window $A(W)$ relative to the production of the certificate for $(T,T')$ one time takes $\frac{1}{2}aa!$ calls for strings of length $\leq\frac{n}{a}$, and we need to apply also some of the routines in Corollary~\ref{co:schreier}, which take time $O(n^{10})$ at most. This can happen at most $n$ times for each window (see the end of \cite[\S 6.1.1]{Hel19}), and the number of windows to update is $\leq m^{2a}$ (see \cite[\S 6.1.2]{Hel19}), so we obtain the claimed runtime for producing the certificates of fullness. Then, we need to generate $F$: we simply take the union of the generators of all certificates, but we do it one certificate at a time and we apply Schreier-Sims at every step, so that the number of generators stays quadratic in $n$ (see the observation after Proposition~\ref{pr:schreier}). The certificates of fullness are at most $m^{2a}$, so this cost is absorbed in the additional time already.

Finally, we need to justify the bounds on $a$ given in the statement. First, by the restrictions on $m,n$ we must have $m,n\geq X$, where $X=\sage{mlim}$ in the CFSG case and $X=\sage{Clim_free}e^{1/\varepsilon^{2}}$ (say) in the CFSG-free case: these are the choices we made in Remark~\ref{re:thism}. The conditions then follow, noticing that for our choice of $n$ the two $a$ respect all bounds in Lemma~\ref{le:epicfsg} and for $\varepsilon$ small the two intervals are large enough to contain an integer.
\end{proof}

Again, $\varepsilon<\frac{1}{100}$ is plenty enough.

Let us also insert here a short lemma that we will use as part of the aggregation of certificates: it is a classical bound on $d$-transitivity for non-giants.

\begin{lemma}\label{le:dtrans}
Let $G\leq\mathrm{Sym}(n)$ be $d$-transitive and $G\neq\mathrm{Sym}(n),\mathrm{Alt}(n)$. Then
\begin{enumerate}[(a)]
\item\label{le:dtranscfsg} $d\leq 5$ (assuming CFSG), or
\item\label{le:dtransfree} $d\leq 3\log n$ (without assuming CFSG).
\end{enumerate}
\end{lemma}

\begin{proof}
See \cite[Thm.~4.11]{Cam99} for the CFSG result and \cite[Satz C]{Wie34} for the CFSG-free result.
\end{proof}

Then we estimate the cost of another major routine, the one represented by the \textit{Design Lemma} and \textit{Split-or-Johnson} (see \cite[\S\S 6-7]{Bab16a} or \cite[\S 5]{Hel19}\footnote{In parts of the next proof, we use terms from the English version \cite{HBD17} instead of the original French ones. The author thinks the reader is better served by this choice, considering also that Babai's original article is in English.}).

\begin{mysage}
fmtrue=3*sqrtinterval(RIF(6))*(3*loginterval(loginterval(RIF(mlim)))+3*loginterval(RIF(6))-2)/sqrtinterval(loginterval(RIF(mlim)))+(3*loginterval(loginterval(RIF(mlim)))+loginterval(6^3*4*RIF(pi)^2))/(4*loginterval(RIF(mlim))^2)+1/(72*sqrtinterval(6))/(loginterval(RIF(mlim)))^(7/2)
fm=roundup(fmtrue,5)
fmsoj=ceil(RIF(fm)+12/loginterval(RIF(3/2)))
\end{mysage}

\begin{proposition}\label{pr:computsoj}
Let $\mathfrak{X}$ be a $b$-ary coherent configuration on $\Gamma$, with $|\Gamma|=m\geq\sage{mlim}$ and $2\leq b\leq\frac{1}{2}m$, such that there is no twin class with $>\frac{1}{2}m$ elements. Then we can find either
\begin{enumerate}[(a)]
\item\label{pr:computsojsplit} a coloured $\frac{2}{3}$-partition of $\Gamma$, or
\item\label{pr:computsojjohn} a Johnson scheme of size $\geq\frac{2}{3}m$ inside $\Gamma$,
\end{enumerate}
at a multiplicative cost of $m^{b+\sage{fmsoj}\log m}$ and at an additive cost of $O(m^{b+14})$.
\end{proposition}

Again, the condition on $m$ is the largeness condition of Remark~\ref{re:thism} (regardless of our position on CFSG).

\begin{proof}
As in Proposition~\ref{pr:computcert}, we are going to discuss only the runtime here. The proof of the rest of the statement is contained in \cite[\S\S 5.1-5.2]{Hel19}. The ``multiplicative cost'' we incur here is the cost of fixing images of a certain number of points of $\Gamma$ (or parts of a partition of $\Gamma$, but fixing the image of a point in the part implies fixing the image of the whole part): arbitrarily fixing a point $x$ in a configuration (or in a graph) in an isomorphism problem translates to trying all possible images of that point, consequently multiplying its contribution. See also Remark~\ref{re:multcost}.

First, we plug the configuration $\mathfrak{X}$ into the Design Lemma, so that we can pull out a classical configuration to use inside Split-or-Johnson: this involves a multiplicative cost of $m^{b-1}$ at most, and a time of $O(m^{b})$ to find the right tuple to use (see \cite[\S 5.1]{Hel19}). Then, either we terminate by fixing $1$ more point (i.e.\ another multiplicative cost of $m$) if the new configuration is not primitive, or we call Split-or-Johnson (SoJ, \cite[Thm.~5.3]{Hel19}).

SoJ itself fixes $1$ element and then, if it does not terminate, calls Bipartite Split-or-Johnson (BSoJ, \cite[Prop.~5.7]{Hel19}). Call $T(m,v)$ the number of elements fixed by BSoJ when $|V_{2}|=v$. The base case is $v\leq(6\log m)^{\frac{3}{2}}$, and here the multiplicative cost is at most $v!$; we use Robbins's bound \cite{Rob55} for factorials,
\begin{equation*}
v!<\sqrt{2\pi}v^{v+\frac{1}{2}}e^{-v+\frac{1}{12v}}
\end{equation*}
(the latter being an increasing function), and the cost is in turn bounded by
\begin{equation*}
\sqrt{2\pi}(6\log m)^{\frac{3}{2}\left((6\log m)^{\frac{3}{2}}+\frac{1}{2}\right)}e^{-(6\log m)^{\frac{3}{2}}+\frac{1}{12}(6\log m)^{-\frac{3}{2}}}=m^{f(m)\log m},
\end{equation*}
where
\begin{equation*}
f(m)=\frac{3\sqrt{6}(3\log\log m+3\log 6-2)}{\sqrt{\log m}}+\frac{3\log\log m+\log(6^{3}\cdot 4\pi^{2})}{4\log^{2}m}+\frac{(72\sqrt{6})^{-1}}{\log^{\frac{7}{2}}m}.
\end{equation*}
Now suppose we are outside the base case; first, we apply the Design Lemma again, for a cost of at most
\begin{equation*}
v^{6\left\lceil\frac{\log m}{\log v}\right\rceil}<v^{12\frac{\log m}{\log v}}=m^{12}.
\end{equation*}
Then we fall again into two subcases: either we recur to a new $v$ that is $\leq\frac{2}{3}$ times the old $v$, with no other cost along the way, or we pass through Coherent Split-or-Johnson (CSoJ, \cite[Prop.~5.8]{Hel19}) and recur to $\leq\frac{1}{2}$ times the old $v$, with $1$ more element fixed in the process (in both cases, it might also happen that we exit the recursion, which is even better). The two situations lead to bounds $T(m,v)\leq m^{12}T\left(m',\frac{2}{3}v\right)$ and $T(m,v)\leq m^{13}T\left(m',\frac{1}{2}v\right)$ respectively, where $m'$ may be smaller than $m$ but still $>\frac{2}{3}m$, or we would exit the recursion again. Since $v<m$ and given the bound in the base case, we obtain in the end
\begin{equation*}
T(m,v)\leq m^{f(m)\log m}\cdot\max\left\{m^{12\log_{3/2}m},m^{13\log_{2}m}\right\}=m^{\left(f(m)+\frac{12}{\log 3/2}\right)\log m}.
\end{equation*}
As for the additive time incurred during the procedure, the heaviest costs come from the use of the Weisfeiler-Leman algorithm inside BSoJ (\cite[Alg.~3]{Hel19}, see also \cite{WL68}), which is performed on a $c$-ary configuration of $V_{2}$ with $c\leq 6\left\lceil\frac{\log m}{\log v}\right\rceil$, entailing spending $O(c^{2}v^{2c+1}\log v)\leq O(m^{13}\log^{3}m)$ time for each encounter we have with Weisfeiler-Leman: by what we described before, we call BSoJ at most $O(\log m)$ times, so that we can safely bound the runtime by $O(m^{14})$. All other costs inside SoJ and its relatives (finding twins, colours, etc...) can also be bounded by $O(m^{14})$.

Hence, at the end we incurred a multiplicative cost of $m^{b+\left(f(m)+\frac{12}{\log 3/2}\right)\log m}$ and an additive cost of $O(m^{b+14})$. For $m\geq\sage{mlim}$ we have $f(m)\leq\sage{fm}$, and we obtain the bound in the statement.
\end{proof}

\section{The algorithm}\label{se:algo}

During the whole process, we are working with a pair of strings of the same length $|\Omega|$ and with a group $G$ that respects a system of blocks in $\Omega$; every time we go through the various steps, we are going to either decrease the length of $\Omega$, increase the size of the blocks or decrease the degree of $G$ (in the sense that $G$ will not vary but we will decrease $m$ where $G\leq\mathrm{Sym}(m)$ as abstract groups).

\begin{remark}\label{re:nsmall}
The case of $n$ small is trivial to examine, and could work as a base case for our algorithm (although we actually follow another path): if $n\leq C$ for some fixed constant $C$, then we can determine $\mathrm{Iso}_{G}(\mathbf{x},\mathbf{y})$ in constant time with constant number of generators. 

To achieve this, just try all the permutations of $G$: we can write all its elements in constant time by Corollary~\ref{co:schreier}, then check whether each of them sends $\mathbf{x}$ to $\mathbf{y}$. If we do not find one, $\mathrm{Iso}_{G}(\mathbf{x},\mathbf{y})$ is empty, otherwise after we find the first one (call it $\tau$) we check which elements of $G$ fix $\mathbf{x}$; the collection of all those that pass the test are all the elements of $\mathrm{Aut}_{G}(\mathbf{x})$, and they also trivially form a set of generators of $\mathrm{Aut}_{G}(\mathbf{x})$: since $\mathrm{Iso}_{G}(\mathbf{x},\mathbf{y})=\mathrm{Aut}_{G}(\mathbf{x})\tau$ by Remark~\ref{re:iso} (or by Lemma~\ref{le:chain}\eqref{le:chain1} and $G\tau=G$), we are done.
\end{remark}

As we already mentioned, the base case of the atoms ($\mathcal{A}$) will be treated in a different way, as presented in Proposition~\ref{pr:base}. Here we need only to cover $n=1$, which is trivial: this is also an atom, as $\mathrm{Sym}(1)=\mathrm{Alt}(1)=\{\mathrm{Id}_{1}\}$; from now on we can suppose $n>1$.

Let us start now with the simplest of recursions, the one with $G$ intransitive.

\begin{proposition}\label{pr:Gintrans}
Let $|\Omega|=n$, $G\leq\mathrm{Sym}(\Omega)$ and let $\mathbf{x},\mathbf{y}:\Omega\rightarrow\Sigma$ be two strings. If $G$ is intransitive, we can reduce the problem of determining $\mathrm{Iso}_{G}(\mathbf{x},\mathbf{y})$ to determining sets $\mathrm{Iso}_{G_{i}}(\mathbf{x}_{i},\mathbf{y}_{i})$ such that $\sum_{i}|\mathbf{x}_{i}|=\sum_{i}|\mathbf{y}_{i}|=n$ and each $G_{i}$ is transitive. The reduction takes time $O(n^{11})$ and no multiplicative cost.
\end{proposition}

\begin{proof}
Let $\Delta$ be an orbit induced by the action of $G$ on $\Omega$, nonempty and properly contained in $\Omega$ since $G$ is intransitive; we can find orbits in time $O(n^{3})$ by Lemma~\ref{le:orbblo}. We call $G_{1}=G|_{\Delta},\mathbf{x}_{1}=\mathbf{x}|_{\Delta},\mathbf{y}_{1}=\mathbf{y}|_{\Delta}$ the restriction of $G,\mathbf{x},\mathbf{y}$ to $\Delta$, as in Lemma~\ref{le:chain}\eqref{le:chain4}; we suppose that we can compute the set $\mathrm{Iso}_{G_{1}}(\mathbf{x}_{1},\mathbf{y}_{1})=H_{1}\tau_{1}$. As in Lemma~\ref{le:chain}\eqref{le:chain4}, we will use $\overline{\alpha}$ to indicate the object (or \textit{an} object) whose restriction to a subset of $\Omega$ is $\alpha$: this subset will be either $\Delta$ or $\Omega\setminus\Delta$, depending on $\alpha$; by Corollary~\ref{co:schreier}\eqref{co:schreierpreim} with $s=h=2$, finding $\overline{\alpha}$ from $\alpha$ takes time $O(n^{10})$.

First, by Lemma~\ref{le:chain}\eqref{le:chain4} we have $\mathrm{Iso}_{G}^{\Delta}(\mathbf{x},\mathbf{y})=\overline{H_{1}}\overline{\tau_{1}}$; then, by Lemma~\ref{le:chain}\eqref{le:chain3},
\begin{equation}\label{eq:intranstorest}
\mathrm{Iso}_{G}(\mathbf{x},\mathbf{y})=\mathrm{Iso}_{\overline{H_{1}}}^{\Omega\setminus\Delta}(\mathbf{x},\mathbf{y}^{\overline{\tau_{1}}^{-1}})\overline{\tau_{1}}.
\end{equation}
If we can compute
\begin{equation}\label{eq:intranstherest}
\mathrm{Iso}_{\overline{H_{1}}|_{\Omega\setminus\Delta}}(\mathbf{x}|_{\Omega\setminus\Delta},\mathbf{y}^{\overline{\tau_{1}}^{-1}}|_{\Omega\setminus\Delta})=K_{1}\upsilon_{1},
\end{equation}
we can use again Lemma~\ref{le:chain}\eqref{le:chain4} to plug \eqref{eq:intranstherest} inside \eqref{eq:intranstorest} and obtain that $\mathrm{Iso}_{G}(\mathbf{x},\mathbf{y})=\overline{K_{1}}\overline{\upsilon_{1}}\overline{\tau_{1}}$. The whole process reduces in time $O(n^{10})$ the determination of $\mathrm{Iso}_{G}(\mathbf{x},\mathbf{y})$ to the determination of $\mathrm{Iso}$ sets on the shorter pieces $\Delta,\Omega\setminus\Delta$.

We can repeat the same procedure on the $\mathrm{Iso}$ in \eqref{eq:intranstherest}: notice that the group and the strings are all defined on $\Omega\setminus\Delta$, so if the group $\overline{H_{1}}|_{\Omega\setminus\Delta}$ is intransitive we again have a $\Delta'\subsetneq\Omega\setminus\Delta$, a group $G_{2}=\overline{H_{1}}|_{\Delta'}$ and strings $\mathbf{x}_{2}=\mathbf{x}|_{\Delta'},\mathbf{y}_{2}=\mathbf{y}^{\overline{\tau_{1}}^{-1}}|_{\Delta'}$ and we continue as before. This happens at most $n$ times.

In the end, we have spent time $O(n^{11})$ and computed sets $\mathrm{Iso}_{G_{i}}(\mathbf{x}_{i},\mathbf{y}_{i})$: each $G_{i}$ is defined in a way that makes it transitive, because we always restrict to an orbit, and each $\mathbf{x}_{i},\mathbf{y}_{i}$ is the restriction of strings $\mathbf{x},\mathbf{y}^{\sigma}$ to a different part of $\Omega$, so that the sum of their lengths is $n$.
\end{proof}

The partition of $\Omega$ into the orbits of the action of $G$, and the reduction of the problem of determining $\mathrm{Iso}_{G}(\mathbf{x},\mathbf{y})$ to problems on shorter strings, corresponds (in reverse, so to speak) to the glueing process of cosets on disjoint sets featured in ($\mathcal{C}$2).

Then, let us continue tackling the next route to recursion, the case of $G$ imprimitive.

\begin{proposition}\label{pr:Gimprim}
Let $|\Omega|=n$, $G\leq\mathrm{Sym}(\Omega)$ and let $\mathbf{x},\mathbf{y}:\Omega\rightarrow\Sigma$ be two strings. If $G$ is transitive but imprimitive, call $N$ the stabilizer of a minimal set of blocks: then we can reduce the problem of determining $\mathrm{Iso}_{G}(\mathbf{x},\mathbf{y})$ to computing the elements of $G/N$ and determining $|G/N|$ sets $\mathrm{Iso}_{N}(\mathbf{x},\mathbf{y}_{i})$ (where $N$ is intransitive). The reduction takes time $O(|G/N|n^{10})$ and no multiplicative cost.
\end{proposition}

\begin{proof}
Let $\{B_{j}\}_{j}$ be a minimal system of blocks for $G$ (it is not a trivial partition since $G$ is imprimitive), which we can retrieve in time $O(n^{4})$ by Lemma~\ref{le:orbblo}. Let $N$ be the stabilizer of this system: by Corollary~\ref{co:schreier}\eqref{co:schreiersetstab}, we can compute it in time $O(n^{10})$.

Write $G=\bigcup_{i}N\sigma_{i}$, where each $\sigma_{i}$ is a representative of a coset of $N$, so that the number of elements $\sigma_{i}$ is $|G/N|$; if we know all the elements of $G/N$, we can determine each $\sigma_{i}$ in time $O(n^{10})$ by Corollary~\ref{co:schreier}\eqref{co:schreierpreim} with $s=h=2$. By Lemma~\ref{le:chain}\eqref{le:chain1}-\ref{le:chain}\eqref{le:chain2},
\begin{equation*}
\mathrm{Iso}_{G}(\mathbf{x},\mathbf{y})=\bigcup_{i}\mathrm{Iso}_{N}(\mathbf{x},\mathbf{y}^{\sigma_{i}^{-1}})\sigma_{i},
\end{equation*}
so we only have to compute the $\mathrm{Iso}_{N}(\mathbf{x},\mathbf{y}^{\sigma_{i}^{-1}})$ now; after having done so, we have a description of those sets as $H\tau_{i}$ where $H=\mathrm{Aut}_{N}(\mathbf{x})$ is generated by a certain set $S$, and
\begin{equation*}
\mathrm{Iso}_{G}(\mathbf{x},\mathbf{y})=\bigcup_{i}H\tau_{i}\sigma_{i}=\langle S\cup\{\tau_{i}\sigma_{i}\sigma_{1}^{-1}\tau_{1}^{-1}\}_{i}\rangle\tau_{1}\sigma_{1}.
\end{equation*}
Finally, we can filter the set $S\cup\{\tau_{i}\sigma_{i}\sigma_{1}^{-1}\tau_{1}^{-1}\}_{i}$ using the Schreier-Sims algorithm in time $O(n^{5}+n^{3}(n^{2}+|G/N|))$ to obtain a description of $\mathrm{Iso}_{G}(\mathbf{x},\mathbf{y})$ with quadratically many generators, and the claim is proved.
\end{proof}

This process, which essentially reduces the problem to a case-by-case examination, corresponds in reverse to the union of cosets featured in ($\mathcal{C}$1). Proposition~\ref{pr:Gimprim} cannot be used directly, as a case-by-case reduction is very expensive in general: nevertheless, seeing this reduction process is useful, as it is used when $G/N$ is especially small (Corollary~\ref{co:cammar}\eqref{co:cammarsmall}, Proposition~\ref{pr:pybersmall}).

Before going to the key steps of the main algorithm, we introduce a couple of combinatorial lemmas that will be useful in the future. The spirit behind them is to be able to start with the set $\binom{\Gamma}{k}$ of all the $k$-subsets of some $\Gamma$ and:
\begin{enumerate}[(a)]
\item in one case, after finding a partition of $\Gamma$, transfer the partition to $\binom{\Gamma}{k}$ itself (Lemma~\ref{le:parts});
\item in the other case, after identifying $\Gamma$ with another $\binom{\Gamma'}{k'}$, use this identification to partition $\binom{\Gamma}{k}$ (Lemma~\ref{le:johnjohn}).
\end{enumerate}

In the following, a \textit{coloured partition} of a set is a partition in which each part is assigned a colour. A permutation subgroup \textit{respects} a coloured partition if it respects both the partition and the colouring: in other words, for any permutation in the group, the image of any part of a given colour is another part of the same colour.

\begin{mysage}
mlim_for_parts=1046
\end{mysage}

\begin{lemma}\label{le:parts}
Let $|\Gamma|=m$ and let $\mathcal{B}=\binom{\Gamma}{k}$, with $k\leq\sqrt{\frac{m}{\log m}}$; suppose that $G\leq\mathrm{Sym}(\Gamma)$ acts on $\Gamma$ in such a way that there is a coloured partition $\mathcal{C}$ of $\Gamma$ respected by $G$ and whose parts are of size $\leq\alpha|\Gamma|$ (for some $\alpha\leq\frac{2}{3}$). Then either $m\leq\sage{mlim_for_parts-1}$ or $\mathcal{B}$ has a coloured partition $\mathcal{C}'$, respected by the natural action of $G$ on $\mathcal{B}$, whose parts are of size $\leq\frac{2}{3}|\mathcal{B}|$.
\end{lemma}

\begin{proof}
Starting from the partition $\mathcal{C}$ of $\Gamma$, we can naturally construct the following partition $\mathcal{C}'$ of $\mathcal{B}$: each part of $\mathcal{C}'$ collects the elements of $\mathcal{B}$ (i.e.\ the $k$-subsets of $\Gamma$) that intersect each part of $\mathcal{C}$ with a specific intersection size; $\mathcal{C}'$ is also naturally a coloured partition: if in a given part $A'\in\mathcal{C}'$ the ordered tuple of intersection sizes with parts $A_{i}\in\mathcal{C}'$ is $(k_{i})_{i}$, we can give to $A'$ the colour given by the ordered tuple of unordered tuples of intersection sizes for all parts of the same colour for every colour of $\mathcal{C}$ (remember, the fact that $G$ respects $\mathcal{C}$ means that different colours will not mix but different parts of the same colour can be sent to each other).

Now we must prove the claim about the size of the parts $A'_{j}\in\mathcal{C}'$. Fix any part $A'_{0}\in\mathcal{C}'$: from what we said above, all the $k$-subsets belonging to $A'_{0}$ are intersecting the parts of $\mathcal{C}$ in the same number of points, so fix a part $A_{0}\in\mathcal{C}$ whose intersection with them is of a certain size $a>0$. The number of $k$-subsets of $\Gamma$ intersecting $A_{0}$ in $a$ points is $\binom{|A_{0}|}{a}\binom{m-|A_{0}|}{k-a}$, so this is an upper bound for $|A'_{0}|$: we just have to prove that this number is at most $\frac{2}{3}\binom{m}{k}$ (for $m$ large enough).

If $k=1$ the task is already accomplished: in this case in fact we also have $a=1$ and then $|A_{0}|\leq\alpha m\leq\frac{2}{3}m$. From now on, $k>1$.

Let us call $|A_{0}|=\beta m$, where $\beta\leq\alpha\leq\frac{2}{3}$. Then
\begin{align*}
\binom{\beta m}{a}\binom{(1-\beta)m}{k-a} = & \ \frac{1}{k!}\binom{k}{a}\beta m(\beta m-1)\ldots(\beta m-a+1) \ \cdot \\
 & \ \cdot(1-\beta)m((1-\beta)m-1)\ldots((1-\beta)m-k+a+1).
\end{align*}
First, since $\beta<1$ we have obviously $\beta m-i\leq\beta(m-i)$ for all $0\leq i<a$. On the other hand, for $0\leq i<k-a$,
\begin{align*}
\frac{(1-\beta)m-i}{(1-\beta)(m-i-a)} & = 1+\frac{a-(i+a)\beta}{(1-\beta)(m-i-a)} \ \leq \ 1+\frac{a}{m-i-a} \\
 & \leq 1+\frac{k}{m-k} \ < \ 1+\frac{2k}{m},
\end{align*}
so that
\begin{align}
\binom{\beta m}{a}\binom{(1-\beta)m}{k-a} < & \ \frac{1}{k!}\binom{k}{a}\beta^{a}m(m-1)\ldots(m-a+1) \ \cdot \nonumber \\
 & \ \cdot(1-\beta)^{k-a}(m-a)\ldots(m-k+1)\left(1+\frac{2k}{m}\right)^{k-a} \nonumber \\
 = & \ \binom{m}{k}\binom{k}{a}\beta^{a}(1-\beta)^{k-a}\left(1+\frac{2k}{m}\right)^{k-a}. \label{eq:haraldok}
\end{align}
The last factor can be easily bounded in the following way:
\begin{align*}
\left(1+\frac{2k}{m}\right)^{k-a} & <\left(1+\frac{2}{\sqrt{m\log m}}\right)^{k}  \\
 & =\left(1+\frac{2}{\sqrt{m\log m}}\right)^{\frac{\sqrt{m\log m}}{2}\cdot\frac{2k}{\sqrt{m\log m}}} \\
 & <e^{\frac{2}{\log m}}.
\end{align*}
Let us treat the rest now. We are going to prove that
\begin{equation}\label{eq:haraldno}
\binom{k}{a}\beta^{a}(1-\beta)^{k-a}\leq\frac{1}{2}.
\end{equation}
First, we start with the case $k\geq 5$ and $2\leq a\leq k-2$, implying that $a\geq 2$, $k-a\geq 2$ with at least one being a strict inequality. We have
\begin{align*}
\frac{(1-\beta)a}{\beta(k-a+1)}+\frac{\beta(k-a)}{(1-\beta)(a+1)} & = \frac{(1-\beta)^{2}a(a+1)+\beta^{2}(k-a)(k-a+1)}{\beta(1-\beta)(a+1)(k-a+1)} \\
 & > \frac{(1-\beta)^{2}a^{2}+\beta^{2}(k-a)^{2}}{\beta(1-\beta)a(k-a)}\cdot\frac{a}{a+1}\frac{k-a}{k-a+1}.
\end{align*}
The first fraction is of the form $\frac{x^{2}+y^{2}}{xy}$, which is equal to $\frac{(x-y)^{2}}{xy}+2\geq 2$; as for the other two, they are both $\geq\frac{2}{3}$ and at least one is $\geq\frac{3}{4}$: therefore the whole product is $\geq 1$. This means that
\begin{align*}
1 & = (\beta+1-\beta)^{k} \ = \ \sum_{a'=0}^{k}\binom{k}{a'}\beta^{a'}(1-\beta)^{k-a'} \\
 & > \sum_{a'\in\{a-1,a,a+1\}}\binom{k}{a'}\beta^{a'}(1-\beta)^{k-a'} \\
 & = \binom{k}{a}\beta^{a}(1-\beta)^{k-a}\left(\frac{(1-\beta)a}{\beta(k-a+1)}+1+\frac{\beta(k-a)}{(1-\beta)(a+1)}\right) \\
 & \geq 2\binom{k}{a}\beta^{a}(1-\beta)^{k-a},
\end{align*}
and \eqref{eq:haraldno} is proved in this case. For $k=4$ and $a=2$,
\begin{equation*}
\frac{2(1-\beta)}{3\beta}+\frac{2\beta}{3(1-\beta)}=\frac{2}{3}\frac{(1-\beta)^{2}+\beta^{2}}{\beta(1-\beta)}\geq\frac{4}{3}>1,
\end{equation*}
and we are done as before. Now, let $a=1$ or $a=k-1$: we can suppose $a=k-1$ by exchanging the role of $\beta$ and $1-\beta$ if necessary (although we cannot use the bound $\beta\leq\frac{2}{3}$ anymore); $k\beta^{k-1}(1-\beta)$ has a maximum in $\beta=1-\frac{1}{k}$, in which it is equal to $\frac{k}{k-1}\left(1-\frac{1}{k}\right)^{k}$. The factor $\left(1-\frac{1}{k}\right)^{k}$ is bounded from above by $\frac{1}{e}$, so for $k\geq 4$ we obtain the bound $<\frac{1}{2}$; for $k=2,3$ we just check directly obtaining $\frac{1}{2},\frac{4}{9}$ respectively. Finally, let $a=k$: then we have just $\beta^{k}$, which is $\leq\beta^{2}\leq\frac{4}{9}$, and \eqref{eq:haraldno} is proved for all cases.

Plugging our results into \eqref{eq:haraldok},
\begin{equation*}
\binom{\beta m}{a}\binom{(1-\beta)m}{k-a}<\binom{m}{k}\frac{1}{2}e^{\frac{2}{\log m}},
\end{equation*}
and for $m\geq\sage{mlim_for_parts}$ we obtain $\frac{1}{2}e^{\frac{2}{\log m}}<\frac{2}{3}$.
\end{proof}

Given our choice of large $m,n$ inside Remark~\ref{re:thism}, Lemma~\ref{le:parts} applies any time we are assuming $m>C\log^{e}n$ for the appropriate $C,e$.

\begin{mysage}
mplim=12
mplim2=16
\end{mysage}

\begin{lemma}\label{le:johnjohn}
Let $\Gamma'$ be a set, let $\Gamma=\binom{\Gamma'}{k'}$ for some $2\leq k'\leq\frac{|\Gamma'|}{2}$, and let $\mathcal{B}=\binom{\Gamma}{k}$ for some $2\leq k\leq\frac{|\Gamma|}{2}$; suppose that $|\Gamma'|=m'\geq\sage{mplim}$. Let any permutation of $\Gamma'$ induce the natural permutations of $\Gamma$ and $\mathcal{B}$; then any $H\leq\mathrm{Sym}(\Gamma')$ divides $\mathcal{B}$ into a system of orbits and blocks such that each part is $\leq\frac{1}{2}|\mathcal{B}|$.
\end{lemma}

\begin{proof}
Let $\Delta$ be any orbit of $\mathcal{B}$ under the action given in the statement. Any element $x\in\Delta$ is a $k$-set of $k'$-sets of elements of $\Gamma'$: since every $x'\in\Delta$ can be sent to $x$ by some permutation induced by some $h\in H$, all the elements of $\Delta$ are constructed respecting the same equalities among the elements of their elements (for example, if there are $a_{1},a_{2}\in x$ with $b_{1},b_{2},b_{3}\in a_{1}\cap a_{2}$, then any $x'$ also has $a'_{1},a'_{2}$ with $b'_{1},b'_{2},b'_{3}\in a'_{1}\cap a'_{2}$, and so on). Every orbit $\Delta$ is therefore contained in the subset $\mathcal{B}_{r}\subseteq\mathcal{B}$ of elements of $\mathcal{B}$ respecting some given set of relations $r$; if we prove that either $\mathcal{B}_{r}$ is of size $\leq\frac{1}{2}|\mathcal{B}|$ or can be divided into blocks with the same property, the same will hold for $\Delta$ and we would be done.

For any $x\in\mathcal{B}_{r}$, let $A(x)\subseteq\Gamma'$ be the set of the elements of all the elements of $x$, with $|A(x)|=a$ ($a$ does not depend on $x$ since it is determined by the relations $r$); we divide $\mathcal{B}_{r}$ into blocks, where each of them collects all the $x$ with the same $A(x)$: these are really blocks, in the sense that the elements of $\mathcal{B}_{r}$ inside them move together under the action of $H$ since this movement depends ultimately on where $A(x)$ is moved inside $\Gamma'$. We have to exclude that the so formed block system is trivial, i.e.\ that either the blocks have size $1$ or that the whole $\mathcal{B}_{r}$ is a block: if we do it, we are done.

Having blocks of size $1$ means that each $x$ already collects all the possible $k'$-subsets of its own $A(x)$, so that $x$ is its own only permutation under $\mathrm{Sym}(A(x))$: this means that $k=\binom{a}{k'}$ and that $\mathcal{B}_{r}$ has $\binom{m'}{a}$ elements, one for each $A(x)$. $\mathcal{B}$ has $\binom{|\Gamma|}{k}$ elements, where $|\Gamma|=\binom{m'}{k'}$, so to prove the statement in this case it is sufficient to prove that
\begin{equation}\label{eq:johnblock1}
\binom{m'}{a}\leq\frac{1}{2}\binom{\binom{m'}{k'}}{\binom{a}{k'}},
\end{equation}
and we would have shown that $\mathcal{B}_{r}$ is small.

Since $k\geq 2$ there are at least two distinct $k'$-subsets of $\Gamma'$ participating in the formation of $A(x)$, so $a>k'$ and then $a\leq\binom{a}{k'}$; we also recall the easy bounds $\left(\frac{x}{y}\right)^{y}\leq\binom{x}{y}\leq\left(\frac{ex}{y}\right)^{y}$. Then, since $m'\geq\sage{mplim}$, $2\leq k'\leq\frac{m'}{2}$, $k\leq\frac{|\Gamma|}{2}$ and $\left(\frac{11}{2e}\right)^{a}>2$, we obtain
\begin{equation}\label{eq:johncalc}
\binom{\binom{m'}{k'}}{\binom{a}{k'}}\geq\binom{\binom{m'}{k'}}{a}\geq\left(\frac{\binom{m'}{k'}}{a}\right)^{a}\geq\left(\frac{\frac{11}{2}m'}{a}\right)^{a}>2\left(\frac{em'}{a}\right)^{a}\geq 2\binom{m'}{a},
\end{equation}
and \eqref{eq:johnblock1} is proved.

Having $\mathcal{B}_{r}$ as a whole block means that all the $x\in\mathcal{B}_{r}$ are coming from the same $A(x)$; as $\mathcal{B}_{r}$ just collects all elements of $\mathcal{B}$ with the same relations, with no other discriminating condition, $A(x)$ must be the whole $\Gamma'$. For each $x\in\mathcal{B}_{r}$ and $\gamma\in\Gamma'$, call $N(\gamma,x)$ the number of elements of $x$ that contain $\gamma$: the multiset $\{N(\gamma,x)|\gamma\in\Gamma'\}$ is independent from $x$, since it is a reflection of the relations of $\mathcal{B}_{r}$.

Suppose first that such multiset has all equal elements, i.e.\ every $\gamma$ is contained in the same number $N$ of $k'$-subsets of $\Gamma'$ belonging to a fixed $x$ (or to any $x$, given our hypotheses): this is a rather constraining condition in $\mathcal{B}$, so we will show that $\mathcal{B}_{r}$ is small. Consider the set $\mathcal{C}_{1}\subseteq\mathcal{B}$ of all $x$ with multiset $\{N,N,N,\ldots,N\}$ ($m'$ times), so that $\mathcal{B}_{r}\subseteq\mathcal{C}_{1}$, and consider the set $\mathcal{C}_{2}\subseteq\mathcal{B}$ of all $x$ with multiset $\{N+1\ldots,N+1,N-1\ldots,N-1,N,\ldots,N\}$, where the number $k''$ of $N+1$ is equal to the number of $N-1$ and runs among all $1\leq k''\leq k'$: construct the bipartite graph $\mathcal{C}_{1}\cup\mathcal{C}_{2}$ where $\{x_{1},x_{2}\}$ is an edge if and only if we can change exactly one $k'$-subset inside $x_{1}$ to obtain $x_{2}$. Every $x_{1}\in\mathcal{C}_{1}$ has $k\left(\binom{m'}{k'}-k\right)\geq\binom{m'}{k'}$ neighbours, since we can move each of the $k'$-subsets of $x_{1}$ to any of the $k'$-subsets that are not already in $x_{1}$ and obtain some (distinct) element of $\mathcal{C}_{2}$; on the other hand, the number of neighbours of a given $x_{2}$ is at most $\binom{m'-2k''}{k'-k''}$: in fact, each $k'$-subset that contains all the $\gamma$ with $N+1$ can be moved only in one way to produce an element of $\mathcal{C}_{1}$, namely by replacing the $\gamma$ with $N+1$ with the $\gamma$ with $N-1$ and fixing the other ones, and the number of such subsets is bounded by $\binom{m'-2k''}{k'-k''}$. Provided that $b_{i}\leq\frac{1}{2}a_{i}$, $a_{1}\leq a_{2}$ and $b_{1}\leq b_{2}$ imply $\binom{a_{1}}{b_{1}}\leq\binom{a_{2}}{b_{2}}$; therefore
\begin{equation*}
\binom{m'}{k'}|\mathcal{C}_{1}|\leq|\{\mathrm{edges \ of \ }\mathcal{C}_{1}\cup\mathcal{C}_{2}\}|\leq\binom{m'-2k''}{k'-k''}|\mathcal{C}_{2}|\leq \binom{m'}{k'}|\mathcal{C}_{2}|,
\end{equation*}
and since $\mathcal{C}_{1}$ and $\mathcal{C}_{2}$ are disjoint we obtain $|\mathcal{B}_{r}|\leq\frac{1}{2}|\mathcal{B}|$.

Now suppose that the multiset $\{N(\gamma,x)|\gamma\in\Gamma'\}$ has at least two distinct elements; take the least frequent of these elements (or the smallest of the least frequent ones, if more than one exists), say that there are $k''$ of them with $k''\leq\frac{m'}{2}<\frac{kk'}{2}$: the second inequality comes from the fact that $A(x)=\Gamma'$, implying that $kk'\geq m'$, and that equality is excluded because it would imply $N(\gamma,x)=1$ regardless of $\gamma$. Call $A'(x)$ the set of $\gamma$ with this specified $N$ for $x$; $A'(x)$ is properly contained in $\Gamma'$, so there must exist elements $x$ with different $A'(x)$: we collect elements $x\in\mathcal{B}$ based on their $A'(x)$, and as we said before for $A(x)$ this forms a system of blocks, which are not the whole $\mathcal{B}$ since $A'(x)\neq\Gamma'$. We have to exclude that this system has blocks of size $1$.

Assume that these blocks have indeed size $1$, which means that $|\mathcal{B}_{r}|=\binom{m'}{k''}$ (one element for each $A'(x)$); as before, we have to prove that
\begin{equation*}
\binom{m'}{k''}\leq\frac{1}{2}\binom{\binom{m'}{k'}}{k}.
\end{equation*}
When $k''\leq k'$ we have $\binom{\binom{m'}{k'}}{k}\geq\binom{\binom{m'}{k''}}{k}>2\binom{m'}{k''}$, and when $k''\leq k$ we can say $\binom{\binom{m'}{k'}}{k}\geq\binom{\binom{m'}{k'}}{k''}$ and continue as in \eqref{eq:johncalc}, so we can assume $k''>k,k'$; this also excludes the cases $k=2$ and $k'=2$, using $k''<\frac{kk'}{2}$. Let us start with the case $\frac{m'}{k'}>4$; using $k\geq\lceil\frac{m'}{k'}\rceil$, the bounds on binomial coefficients and $m'\geq\sage{mplim}$,
\begin{equation*}
\binom{\binom{m'}{k'}}{k}\geq\binom{\binom{m'}{k'}}{\lceil\frac{m'}{k'}\rceil}\geq\left(\frac{m'}{k'}\right)^{\frac{(k'-1)m'}{k'}}>4^{\frac{2}{3}m'}>2(\sqrt{2e})^{m'}\geq 2\binom{m'}{\lfloor\frac{m'}{2}\rfloor}\geq 2\binom{m'}{k''}.
\end{equation*}
Similarly, for $3<\frac{m'}{k'}\leq 4$ (implying $k\geq 4$),
\begin{equation*}
\binom{\binom{m'}{k'}}{k}\geq\binom{\binom{m'}{k'}}{4}\geq 4^{m'-4}\geq 4^{\frac{2}{3}m'}>2\binom{m'}{k''}.
\end{equation*}
For $2\leq\frac{m'}{k'}\leq 3$ and $k\geq 4$,
\begin{equation*}
\binom{\binom{m'}{k'}}{k}\geq\binom{\binom{m'}{k'}}{4}\geq \frac{3^{\frac{4}{3}m'}}{4^{4}}\geq 4^{\frac{2}{3}m'}>2\binom{m'}{k''}.
\end{equation*}
Finally, for $2\leq\frac{m'}{k'}\leq 3$ and $k=3$, we can first check directly that
\begin{equation*}
\binom{\binom{m'}{k'}}{3}\geq\binom{\binom{m'}{\lceil\frac{m'}{3}\rceil}}{3}\geq 2\binom{m'}{\lfloor\frac{m'}{2}\rfloor}\geq 2\binom{m'}{k''}
\end{equation*}
for each $\sage{mplim}\leq m'<\sage{mplim2}$, while for $m'\geq\sage{mplim2}$
\begin{equation*}
\binom{\binom{m'}{k'}}{3}\geq 3^{m'-3}>2(\sqrt{2e})^{m'}\geq 2\binom{m'}{k''}.
\end{equation*}
Since $\frac{m'}{k'}\geq 2$ is always true, this covers all cases and concludes the proof.
\end{proof}

We are now at a point where we must introduce the cornerstone of the algorithm, the group-theoretic result thanks to which the branching into different cases starts and the recursion is performed. Actually, as anticipated, we have two of them: Theorem~\ref{th:cammar} assumes CFSG and Theorem~\ref{th:pyber} does not; consequently, henceforth we split our reasoning into two different parts, according to our attitude towards CFSG: the two approaches present many points of contact with each other nonetheless, enough to make the proof of the main theorem virtually the same both times.

\subsection{The algorithm, assuming CFSG}

Let us start immediately with our theoretic main tool.

\begin{mysage}
Mathieu=[MathieuGroup(11).order(),MathieuGroup(12).order(),MathieuGroup(23).order(),MathieuGroup(24).order()]
C0=max(Mathieu)
\end{mysage}

\begin{theorem}\label{th:cammar}
Let $|A|=a$ and let $G\leq\mathrm{Sym}(A)$. Assume CFSG. If $G$ is primitive, then one of the following alternatives holds:
\begin{enumerate}[(a)]
\item\label{th:cammarsmall} $|G|\leq C(a)=\max\{C_{0},a^{1+\log_{2}a}\}$ for $C_{0}=\sage{C0}$;
\item\label{th:cammaralt} there is a system $\mathcal{A}$ of (possibly size $1$) blocks of $A$ with $|\mathcal{A}|=\binom{b}{t}\leq a$ and there is a $G'\unlhd G$ with $[G:G']\leq a$ and preserving $\mathcal{A}$, such that we can construct in time $O(n^{10})$ a bijection $\varphi$ between $\mathcal{A}$ and the set $\binom{B}{t}$ of $t$-subsets of a $b$-set $B$ in a way that makes $G'$ isomorphic to $\mathrm{Alt}(B)$, with the action of $G'$ on $\mathcal{A}$ agreeing with the natural action induced by $\mathrm{Alt}(B)$ on $\binom{B}{t}$.
\end{enumerate}
\end{theorem}

\begin{proof}
This theorem is a consequence of Cameron's classification of primitive permutation groups in its formulation due to Mar\'oti \cite{Mar02}. Case \eqref{th:cammarsmall} in the present result collects cases (ii) and (iii) in \cite[Thm.\ 1.1]{Mar02}, and $C_{0}$ is the size of the largest of the four Mathieu groups that appear in (ii), namely $\mathrm{M}_{24}$. The other alternative is realized by taking case (i) and choosing $G'$ to be the power of $\mathrm{Alt}(B)$ that is guaranteed to exist as a subgroup of $G$; the rest of the structure is retrieved by creating the partition $\mathcal{A}$ with one block for each of the possible values of the first coordinate (say) in the formulation of the wreath product in Definition~\ref{de:wreath}, and forgetting the structure coming from all other coordinates, so that we see only one $\mathrm{Alt}(B)$ among all the ones that compose $G'$.

As for the polynomial-time construction of $\varphi$, it is described in \cite[\S 4]{BLS87} (see also \cite[\S 2.8]{Hel19}). The procedure NATURAL\_ACTION thereby described produces a set $D$ divided into blocks $\{B_{i}|i\in I\}$ such that the elements of $A$ correspond to subsets of $D$ of a certain form; our $B$ is any of the $B_{i}$ (say $B_{1}$) and if $\pi:G\rightarrow\mathrm{Sym}(I)$ is the map describing how $G$ permutes the $B_{i}$ then our $G'$ is $\pi^{-1}(\mathrm{Sym}(I)_{(1)})$. All the passages involved in finding $B$ and $G'$ and constructing $\varphi$ come from Corollary~\ref{co:schreier} and Lemma~\ref{le:orbblo} (on sets of size at most $|A|^{2}$): together, they cost at most time $O(n^{10})$ as claimed.
\end{proof}

When we start the whole algorithm to compute $\mathrm{Iso}_{G}(\mathbf{x},\mathbf{y})$, we can divide $G$ into its orbits and blocks (if $G$ is intransitive or imprimitive) in time $O(n^{4})$ by Lemma~\ref{le:orbblo}, and then treat the intransitive case thanks to Proposition~\ref{pr:Gintrans}: therefore we can suppose that $G$ is transitive and acts primitively on some system of blocks $\mathcal{B}$ that we are able to assume to be known.

\begin{corollary}\label{co:cammar}
Let $|\Omega|=n$, $G\leq\mathrm{Sym}(\Omega)$ and let $\mathbf{x},\mathbf{y}:\Omega\rightarrow\Sigma$ be two strings; let $\mathcal{B}$ be a system of blocks of $\Omega$ with $1<|\mathcal{B}|=r\leq n$, on which $G$ acts primitively: call $N$ the stabilizer of the system $\mathcal{B}$, and suppose that there are a set $\Gamma$ of size $m$ and a bijection between $\mathcal{B}$ and $\binom{\Gamma}{k}$ (for some $k$) such that the action of $G/N$ on $\mathcal{B}$ corresponds to the action of some transitive subgroup $H\leq\mathrm{Sym}(\Gamma)$ on $\binom{\Gamma}{k}$. Assume CFSG. Then we can reduce the problem of determining $\mathrm{Iso}_{G}(\mathbf{x},\mathbf{y})$ to one of the following problems:
\begin{enumerate}[(a)]
\item\label{co:cammarsmall} determining $\leq m^{\sage{Clim}\log^{2}n}$ sets of isomorphisms $\mathrm{Iso}_{M}(\mathbf{x},\mathbf{y}_{i})$, where $M\unlhd N$ stabilizes all blocks, in time $O(m^{\sage{Clim}\log^{2}n}n^{10})$ and at no multiplicative cost;
\item\label{co:cammarlem} determining $\leq m$ sets of isomorphisms $\mathrm{Iso}_{G'}(\mathbf{x},\mathbf{y}_{i})$, where $G'$ respects a system of orbits and/or blocks $\mathcal{B}'$ strictly coarser than $\mathcal{B}$ and whose parts are of size $\leq\frac{2}{3}|\Omega|$, in time $O(n^{10})$ and at no multiplicative cost;
\item\label{co:cammaralt} determining $\leq m$ sets of isomorphisms $\mathrm{Iso}_{G'}(\mathbf{x},\mathbf{y}_{i})$, where $G'/N$ acts on $\mathcal{B}$ in the same way as $\mathrm{Alt}(\Gamma')$ acts on $\binom{\Gamma'}{k'}$ (where $|\Gamma'|=m'>\sage{Clim}\log^{2}n$), in time $O(n^{10})$ and at no multiplicative cost.
\end{enumerate}
\end{corollary}

\begin{proof}
Before we start, we point out that we hypothesize the existence of $\Gamma$ in the statement (or, from another perspective, the fact that $k$ may be $\geq 2$) because we want to leave open the possibility that we are returning to this situation after having already been through this step before and found a bijection as in Theorem~\ref{th:cammar}\eqref{th:cammaralt} (using the theorem itself or by other means) that we have then carried forth until this moment, as it may happen. In any case, either we are provided with such $\Gamma,k,\mathcal{B},N$ from past procedures, or in their absence we can determine $\mathcal{B},N$ in time $O(n^{10})$ by Lemma~\ref{le:orbblo} and Corollary~\ref{co:schreier}\eqref{co:schreiersetstab} (setting $\mathcal{B}=\{\{x\}|x\in\Omega\}$ if $G$ is primitive) and then impose $\Gamma=\mathcal{B}$ and $k=1$.

As it can be imagined, we want to use Theorem~\ref{th:cammar} on $A=\Gamma$. First, $H$ must be primitive: if it were not, then its action on $\binom{\Gamma}{k}$ would also be imprimitive (even intransitive, if $k>1$) and this contradicts our hypothesis on $G$; hence we can actually use the theorem. The generators of $G$ (at most $n^{2}$ in number) can be seen as generators of $G/N\simeq H$ and can be processed through Schreier-Sims to determine $|H|$ in time $O(n^{5})$ by Corollary~\ref{co:schreier}\eqref{co:schreiersize}, so that we are able to determine whether we are in case \eqref{th:cammarsmall} or \eqref{th:cammaralt} of Theorem~\ref{th:cammar}.

If we are in case \eqref{th:cammarsmall}, we can write all the elements of $H$ in time $O(n^{5}+C(m)n^{2})$ by Corollary~\ref{co:schreier} and we are exactly in the situation described in Proposition~\ref{pr:Gimprim} (with the computation of all the elements of $H\simeq G/N$ already taken care of). This falls into case \eqref{co:cammarsmall} of the present corollary: we have $N=M$ for the subgroup; also, for $n\leq 3$ obviously $|G/N|\leq m!\leq m^{\sage{Clim}\log^{2}n}$, while for $n\geq 4$ both $C_{0}<2^{\sage{Clim}\log^{2}4}\leq m^{\sage{Clim}\log^{2}n}$ and $m^{1+\log_{2}m}<m^{\sage{Clim}\log^{2}n}$, so the bound on the number of problems holds. The runtime, in light of the previous reasoning on $C(m)$, is also $O(m^{\sage{Clim}\log^{2}n}n^{10})$ as required.

If we are in case \eqref{th:cammaralt}, there is some $H'\unlhd H$ with $[H:H']\leq m$ acting on a partition $\Gamma^{o}$ of $\Gamma$ as $\mathrm{Alt}(\Gamma')$ acts on $\binom{\Gamma'}{k'}$ for some $|\Gamma'|=m'$ and some $k'\geq 1$: $\Gamma^{o},\Gamma',k'$ and the action are all found in time $O(n^{10})$, as we already said. First, suppose that $m\leq \sage{Clim}\log^{2}n$: then $|G/N|<m^{m}\leq m^{\sage{Clim}\log^{2}n}$, and repeating what we did before we retrieve again case \eqref{co:cammarsmall}.

Now suppose that $m>\sage{Clim}\log^{2}n$ and that $\Gamma^{o}$ is a nontrivial partition: as observed in Remark~\ref{re:thism} we have $m\geq\sage{mlim}$, and the hypothesis on $\Gamma^{o}$ makes it into a coloured partition (with only one colour) whose parts are of size $\leq\frac{1}{2}|\Gamma|$; to use Lemma~\ref{le:parts}, we still have to prove that $k\leq\sqrt{\frac{m}{\log m}}$. For $k=1$ this is true for any $m$, so suppose that $k>1$. Obviously we can assume that $m\geq 2k$: in fact there is a natural identification between $\binom{\Gamma}{k}$ and $\binom{\Gamma}{|\Gamma|-k}$, just by taking the complement of each of their elements; therefore
\begin{equation*}
n\geq\binom{m}{k}\geq\left(\frac{m}{k}\right)^{k}\geq 2^{k} \ \Longrightarrow \ k\leq\frac{1}{\log 2}\log n \ \Longrightarrow \ m>k^{2}\cdot\sage{Clim}\log^{2}2>k^{2},
\end{equation*}
and, using this new bound again,
\begin{equation*}
n\geq\binom{m}{k}\geq\left(\frac{m}{k}\right)^{k}>k^{k} \ \Longrightarrow \ \log n>k\log k.
\end{equation*}
The function $f(y)=\frac{y}{\sqrt{\log y}}$ is increasing and $f(k\log k)>k$ for $k>1$, therefore using $k\log k<\log n\leq\sqrt{\frac{1}{\sage{Clim}}m}$ we get $k<\sqrt{\frac{m}{\sage{Clim}\log\sqrt{\frac{1}{\sage{Clim}}m}}}<\sqrt{\frac{m}{\log m}}$ (where $m\geq\sage{mlim}$ is amply sufficient to satisfy the second inequality). Now we are free to use Lemma~\ref{le:parts}, which makes us fall into case \eqref{co:cammarlem} of the present corollary.

Finally, let us have $m>\sage{Clim}\log^{2}n$ and $\Gamma^{o}=\Gamma$: since $m\geq\sage{mlim}$ and $m=\binom{m'}{k'}$, we have $m'\geq 12$ regardless of our choice of $k'$. If both $k$ and $k'$ are $>1$, we can use Lemma~\ref{le:johnjohn} and we fall again into case \eqref{co:cammarlem}. If $k'=1$, then $\Gamma'=\Gamma$ and $H'$ acts as $\mathrm{Alt}(\Gamma)$ on $\Gamma$ itself, thus acting as $\mathrm{Alt}(\Gamma)$ on $\binom{\Gamma}{k}\simeq\mathcal{B}$. If $k=1$, then $\Gamma=\mathcal{B}$ and $H'$ acts as $\mathrm{Alt}(\Gamma')$ on $\binom{\Gamma'}{k'}\simeq\mathcal{B}$; if $m'\leq\sage{Clim}\log^{2}n$ we reduce again to case \eqref{co:cammarsmall} exactly as before, so $m'>\sage{Clim}\log^{2}n$. In both cases, whether $k'=1$ or $k=1$, we can take the pullback $G'$ of $H'$ in $G$ (in time $O(n^{10})$ by Corollary~\ref{co:schreier}\eqref{co:schreierpreim}) and $G'/N\simeq H'$ will satisfy the requirements of case \eqref{co:cammaralt} of this corollary: in fact $[G:G']=[H:H']$ and we can obtain (a preimage of) all the elements of $G/G'$ in time $O(n^{10})$, continuing then with $\mathrm{Iso}_{G}(\mathbf{x},\mathbf{y})=\bigcup_{i}\mathrm{Iso}_{G'}(\mathbf{x},\mathbf{y}^{\sigma_{i}^{-1}})\sigma_{i}$ as in Proposition~\ref{pr:Gimprim}.
\end{proof}

We point out that \cite{Hel19} uses actually a bound on $m$ of the form $m>C\log n$ for the case equivalent to our case \eqref{co:cammaralt}. In order to follow our line of thought we need a stronger bound, quadratic in $\log n$, because otherwise we obtain a weaker inequality than $k\leq\sqrt\frac{m}{\log m}$ and then Lemma~\ref{le:parts} does not work: the issue is with the last factor in \eqref{eq:haraldok}, which needs to decrease with the growth of $m$; the problem is treated incorrectly in \cite[\S 4.2]{Hel19}. A bound $m>C\log n$ is more than we need to obtain the bound on the runtime of the form $n^{O(\log^{2}n)}$ anyway: as observed in \cite[\S 3.1]{Hel19}, it is consistent even with a $n^{O(\log n)}$ runtime, to this day unproven.

After we have reached case \eqref{co:cammarsmall} in the previous corollary, we can simply go through Proposition~\ref{pr:Gintrans} and reduce to examinate each block singularly: this makes $n$ decrease, and we return to the top of this corollary. After case \eqref{co:cammarlem}, $\Omega$ is divided into orbits and blocks that are coarser than the original $\mathcal{B}$: this makes $n$ decrease or the block size increase (or both). Case \eqref{co:cammaralt} is the one we will examine in the following results.

\begin{proposition}\label{pr:base}
Let $|\Omega|=n$, and let the action of $G\leq\mathrm{Sym}(\Omega)$ on $\Omega$ be as in Corollary~\ref{co:cammar}\eqref{co:cammaralt}, i.e.\ there is a system of blocks $\mathcal{B}$ with stabilizer $N$ such that $G/N$ acts on it as $\mathrm{Alt}(\Gamma)$ acts on $\binom{\Gamma}{k}$, where $|\Gamma|=m$ and $|\mathcal{B}|=\binom{m}{k}$. If $k=1$ and the blocks have size $1$, the set $\mathrm{Iso}_{G}(\mathbf{x},\mathbf{y})$ can be determined in time $O(n^{6})$ with at most $n^{2}$ generators.
\end{proposition}

\begin{proof}
Having $k=1$ means that $\Gamma=\Omega$, and having block size $1$ means that $G=G/N\simeq\mathrm{Alt}(\Gamma)=\mathrm{Alt}(\Omega)$. This is a trivial case: if $\mathbf{x}$ and $\mathbf{y}$ do not send the same number of elements of $\Omega$ to the same letter of the alphabet $\Sigma$, the set is empty.

Otherwise, we first obtain $\mathrm{Aut}_{\mathrm{Sym}(\Omega)}(\mathbf{x})$ as a product $\prod_{i}\mathrm{Sym}(\Delta_{i})$, where the $\Delta_{i}$ are the parts of $\Omega$ whose elements are sent by $\mathbf{x}$ to the same letter: more precisely, for each generator of $\mathrm{Sym}(\Delta_{i})$ we find the corresponding element in $\mathrm{Sym}(\Omega)_{(\Omega\setminus\Delta_{i})}$, and then we take the union of these preimages for all $i$; each $\mathrm{Sym}(\Delta_{i})$ can be described by two generators, a transposition and a cycle of length $|\Delta_{i}|$, therefore up until now we are working with $\leq\frac{2}{3}n$ generators. Then, we find $H=\mathrm{Aut}_{\mathrm{Alt}(\Omega)}(\mathbf{x})$: by Corollary~\ref{co:schreier}\eqref{co:schreiersubtest}, since the index is $\leq 2$ and the test to prove whether a permutation is even is linear-time (just by computing the length of the cycles), we obtain polynomially many generators of $H$ in time $O(n^{5})$; more precisely, the number of generators is at most $\left(\frac{2}{3}n+1\right)^{3}$ by Schreier's lemma (\cite{Sch27}, see for example \cite[Lemma~4.2.1]{Ser03}) and we can reduce it to $\leq n^{2}$ using Schreier-Sims and spending time $O(n^{6})$ by Proposition~\ref{pr:schreier}.

Finally we take any bijection $\pi:\Omega\rightarrow\Omega$ sending elements sent to each letter of $\Sigma$ by $\mathbf{x}$ to the elements sent to the same letter by $\mathbf{y}$. If this bijection is in $\mathrm{Alt}(\Omega)$ we have $\mathrm{Iso}_{G}(\mathbf{x},\mathbf{y})=H\pi$; if it is not, there are two possibilities: if there is a letter that appears twice in the strings (say $\mathbf{x}(r_{1})=\mathbf{x}(r_{2})$) we have $\mathrm{Iso}_{G}(\mathbf{x},\mathbf{y})=H\tau\pi$ where $\tau$ is the transposition $(r_{1} \ r_{2})$, otherwise the set is empty again.
\end{proof}

The situation described in Proposition~\ref{pr:base} (apart from the case of $n=1$ taken care of in Remark~\ref{re:nsmall}) is the only true base case of the whole algorithm; the rest of the time, the procedure either stops and gives $\emptyset$ as a result or it reduces to simpler cases, until we arrive to the one given above. Proposition~\ref{pr:base} corresponds to the case of the atom ($\mathcal{A}$) in the main theorem.

Let us see what happens aside from the base case.

\begin{mysage}
multe_cfsg=ceil(fmsoj+alim)
\end{mysage}

\begin{theorem}\label{th:cert}
Let $|\Omega|=n$, let $\mathbf{x}:\Omega\rightarrow\Sigma$ be a string, and let the action of $G\leq\mathrm{Sym}(\Omega)$ on $\Omega$ be as in Corollary~\ref{co:cammar}\eqref{co:cammaralt}, i.e.\ there is a system of blocks $\mathcal{B}$ such that $G$ acts on it as $\mathrm{Alt}(\Gamma)$ acts on $\binom{\Gamma}{k}$, where $|\Gamma|=m$ and $|\mathcal{B}|=\binom{m}{k}$. Assume CFSG; suppose also that $m>\sage{Clim}\log^{2}n$. Then we can reduce to one of the following cases:
\begin{enumerate}[(a)]
\item\label{th:certnodom} $\Gamma$ has a canonical coloured partition in which each part has size $\leq\frac{1}{2}|\Gamma|$;
\item\label{th:certlargesym} there is a canonical set $S\subseteq\Gamma$ of size $>\frac{1}{2}|\Gamma|$ such that for any $\sigma\in\mathrm{Alt}(S)$ there is an element of $\mathrm{Aut}_{G}(\mathbf{x})$ that induces $\sigma$ on $S$;
\item\label{th:certsoj} at a multiplicative cost of at most $m^{\sage{multe_cfsg}\log n}$, either:
\begin{enumerate}[(\ref{th:certsoj}1)]
\item\label{th:certsplit} $\Gamma$ has a coloured partition in which each part has size $\leq\frac{2}{3}|\Gamma|$, or
\item\label{th:certjohn} there are two disjoint sets $V_{1},V_{2}\subseteq\Gamma$, with $V_{2}$ divided into a system of blocks $\mathcal{G}$ with $\binom{|\mathcal{G}|}{k'}=|V_{1}|\geq\frac{2}{3}|\Gamma|$ for some $k'\geq 2$, and there is a bijection between $V_{1}$ and $\binom{\mathcal{G}}{k'}$ such that if a $g\in G$ induces a permutation $\sigma\in\mathrm{Sym}(\mathcal{G})$ of the blocks then it also induces the corresponding permutation of $V_{1}$ through the identification of its elements with the $k'$-subsets of $\mathcal{G}$.
\end{enumerate}
\end{enumerate}
The time necessary for this reduction is the cost of $\frac{1}{2}m^{2a}naa!$ calls of the whole algorithm for strings of length $\leq\frac{n}{a}$ where $a\in\left(\sage{alim_low},\sage{alim}\right)\cdot\log n$, plus some additional time $O(m^{3a}n^{11})$.
\end{theorem}

\begin{proof}
We are in the scenario of Proposition~\ref{pr:computcert}: the construction in our hypothesis yields in particular a surjective map $\phi:G\rightarrow\mathrm{Alt}(\Gamma)$. After $\frac{1}{2}m^{2a}naa!$ calls of the algorithm for strings of length $\leq\frac{n}{a}$ and an additional time of $O(m^{2a}n^{11})$, we have obtained the group $F\leq\mathrm{Aut}_{G}(\mathbf{x})$ generated by all certificates of fullness. Now we follow the case subdivision in \cite[\S 6.2]{Hel19}.
\begin{itemize}
\item ``Cas 1'' is the case of $|\mathrm{supp}(\phi(F))|\geq\frac{1}{2}m$ and no orbit of $\phi(F)$ of length $>\frac{1}{2}m$.
\item ``Cas 2a'' is the case of $|\mathrm{supp}(\phi(F))|\geq\frac{1}{2}m$, an orbit $\Phi$ of $\phi(F)$ of length $>\frac{1}{2}m$, and $\mathrm{Alt}(\Phi)\leq\phi(F)|_{\Phi}$.
\item ``Cas 2b'' is the case of $|\mathrm{supp}(\phi(F))|\geq\frac{1}{2}m$, an orbit $\Phi$ of $\phi(F)$ of length $>\frac{1}{2}m$, and $\mathrm{Alt}(\Phi)\not\leq\phi(F)|_{\Phi}$.
\item ``Cas 3'' is the case of $|\mathrm{supp}(\phi(F))|<\frac{1}{2}m$.
\end{itemize}

In ``Cas 1'' we colour each element of $\Gamma$ by the length of its orbit (in time $O(m^{3})$ by Lemma~\ref{le:orbblo}) and we are in our case \eqref{th:certnodom}. ``Cas 2a'' is our case \eqref{th:certlargesym} for $S=\Phi$.

``Cas 2b'' starts by arbitrarily fixing some points of $\Gamma$, precisely $d-1$ many for $d$ as in Lemma~\ref{le:dtrans}\eqref{le:dtranscfsg}, and then feeds the resulting configuration to the Split-or-Johnson procedure (without passing through the Design Lemma). In ``Cas 3'', the information we already have at hand after the production of the local certificates lets us have a colouring of $(\Gamma\setminus\mathrm{supp}(\phi(F)))^{a}$ with less than half twins (as long as $a\leq\frac{m}{4}$): we can make it into an $a$-ary configuration and refine it through Weisfeiler-Leman at a cost of $O(a^{2}m^{2a+1}\log m)$ for the runtime, and then invoke the Design Lemma plus Split-or-Johnson.

In both cases, we can apply Proposition~\ref{pr:computsoj}: the two alternatives \eqref{pr:computsojsplit} and \eqref{pr:computsojjohn} therein correspond respectively to cases \eqref{th:certsplit} and \eqref{th:certjohn} here. We have explicitly written in our statement what the sentence ``nous pouvons trouver [...] un sch\'ema de Johnson plong\'e sur [...] $\Gamma$'' means in the statement of \cite[Thm.~5.3]{Hel19}: in particular, the fact that the objects that when permuting induce a permutation of $V_{1}$ may be the parts of $\mathcal{B}'$ (instead of being directly the elements of $V_{2}$) is due to the use of \cite[Ex.~2.18]{Hel19} inside CSoJ, where from a graph made of elements of $V_{2}$ we pass to a contracted graph made of its parts.

The multiplicative cost of ``Cas 2b'' and ``Cas 3'' is bounded by $m^{a+\sage{fmsoj}\log m}$ (certainly $d\leq a$, so the ``Cas 2b'' expense is subsumed by the ``Cas 3'' expense), and their additive cost is safely absorbed into the $O(m^{3a}n^{11})$. For our choice of $a$, we obtain the cost featured in \eqref{th:certsoj}.
\end{proof}

\begin{remark}\label{re:multcost}
The multiplicative cost described in case \eqref{th:certsoj} of Theorem~\ref{th:cert} means the following: since a permutation in $G$ induces also an even permutation of $\Gamma$, for any choice of $s$ points $x_{1},\ldots,x_{s}\in\Gamma$ each isomorphism from $\mathbf{x}$ to $\mathbf{y}$ falls into a particular coset of the stabilizer of these points; these cosets are one for each possible choice of images of the points in $\Gamma$.

Call $N$ the preimage in $G$ of $\mathrm{Alt}(\Gamma)_{(x_{1},\ldots,x_{s})}$, found in time $O(n^{5})$ by Corollary~\ref{co:schreier}\eqref{co:schreierpwstab} ($N$ need not be normal in $G$: we call it $N$ in analogy to Proposition~\ref{pr:Gimprim}); $[G:N]\leq m^{s}$, so again by Corollary~\ref{co:schreier}\eqref{co:schreiersubtest} we can write an element $\sigma_{i}$ of each coset of $N$ in time $O(n^{5}+m^{s}n^{3})$. Thus the problem of determining $\mathrm{Iso}_{G}(\mathbf{x},\mathbf{y})$ reduces to $\leq m^{s'}$ problems of determining $\mathrm{Iso}_{N}(\mathbf{x},\mathbf{y}_{i})$, because
\begin{equation*}
\mathrm{Iso}_{G}(\mathbf{x},\mathbf{y})=\bigcup_{i}\mathrm{Iso}_{N}(\mathbf{x},\mathbf{y}^{\sigma_{i}^{-1}})\sigma_{i}
\end{equation*}
exactly as in Proposition~\ref{pr:Gimprim}. It is important to consider that $s'$ as above, the exponent of the multiplicative cost, is not the same as $s$ (despite them being certainly related) and is indeed smaller: the fact is that the elements of $\Gamma$ are not all indistinguishable (due to the presence of $V_{1},V_{2}$), so many possibilities for the choice of $x_{1},\ldots,x_{s}$ are as a matter of fact forbidden; seen in a different light, many of the $\mathrm{Iso}_{N}$ that emerge are known to be empty without the need for computing them, as they do not make $V_{1},V_{2}$ correspond in $\mathbf{x}$ and $\mathbf{y}^{\sigma_{i}^{-1}}$.
\end{remark}

Now that the situation described in the hypothesis of Theorem~\ref{th:cert} has been split into its various cases, we show how to treat each of them while making at least one among our parameters $n,|\mathcal{B}|,m$ decrease.

\begin{corollary}\label{co:certpart}
Let $|\Omega|=n$, $G\leq\mathrm{Sym}(\Omega)$ and let $\mathbf{x},\mathbf{y}:\Omega\rightarrow\Sigma$ be two strings; let $\mathcal{B}$ be a system of blocks such that $G$ acts on it as $\mathrm{Alt}(\Gamma)$ acts on $\binom{\Gamma}{k}$, where $|\Gamma|=m$ and $|\mathcal{B}|=\binom{m}{k}$. Suppose that $m>\sage{Clim}\log^{2}n$; suppose also that, fixing the images $(y_{i})_{i=1}^{s}$ of some elements $(x_{i})_{i=1}^{s}\subseteq\Gamma$, we can find a coloured partition of $\Gamma$ in which each part has size $\leq\alpha|\Gamma|$ (with $\alpha\leq\frac{2}{3}$).

Then, if $N$ is the preimage of $\mathrm{Alt}(\Gamma)_{(x_{1},\ldots,x_{s})}$ inside $G$, $N$ divides $\Omega$ into a system $\mathcal{B}'$ of orbits and blocks (at least as coarse as $\mathcal{B}$) of size $\leq\frac{2}{3}|\Omega|$. Moreover, for any orbit $\Delta$ with $|\Delta|>\frac{2}{3}|\Omega|$, $\mathcal{B}'|_{\Delta}$ is nontrivial and strictly coarser than $\mathcal{B}|_{\Delta}$ and its elements are $k$-subsets of blocks of $\mathcal{B}$ all contained in the same colour $\Gamma_{0}$ of $\Gamma$ of size $>\frac{2}{3}|\Gamma|$; also, the stabilizer of blocks of $\mathcal{B}'|_{\Delta}$ coincides with the stabilizer of blocks of $\Gamma_{0}$.
\end{corollary}

\begin{proof}
This corollary covers cases \eqref{th:certnodom} and \eqref{th:certsplit} of Theorem~\ref{th:cert}. The focus on $N$ is due to the reduction to the problem of determining $\mathrm{Iso}_{N}(\mathbf{x},\mathbf{y}^{\sigma^{-1}})$ featured in Remark~\ref{re:multcost}, where $\sigma\in G$ is an element that sends each $x_{i}$ to $y_{i}$.

We have a coloured partition $\mathcal{C}$ on $\Gamma$ with parts of size $\leq\alpha|\Gamma|$ (with $\alpha\leq\frac{2}{3}$); we can repeat the same reasoning as in Corollary~\ref{co:cammar} (the case $m>\sage{Clim}\log^{2}n$ and $\Gamma^{o}$ nontrivial) and show that the hypotheses of Lemma~\ref{le:parts} hold here. By this lemma, $\Omega$ itself has a coloured partition $\mathcal{C}'$ that is at least as coarse as $\mathcal{B}$ and whose parts are also of size $\leq\frac{2}{3}|\Omega|$: the fact that $N$ respects the colours of $\mathcal{C}'$ means that elements with different colours will not be sent to each other, i.e.\ they sit in different orbits, while respecting the parts with the same colours translates to sending all the elements of one part to the same part, i.e.\ moving them as a block.

If we are in an orbit $\Delta$ of size $>\frac{2}{3}|\Omega|$, it means that inside $\mathcal{C}'$ we are in a colour of size $>\frac{2}{3}|\Omega|$, so that it will also have to be divided into smaller parts with the same colour: therefore, $\mathcal{B}'|_{\Delta}$ is nontrivial and strictly coarser than $\mathcal{B}|_{\Delta}$, since each part will contain not all blocks and at least two blocks of $\mathcal{B}$. Using the reasoning in Lemma~\ref{le:parts}, $\Delta$ must come from a $\Gamma_{0}$ as in our statement, and by our description of $\mathcal{C}'$ in that lemma the block stabilizer of $\mathcal{B}'|_{\Delta}$ contains the block stabilizer of $\Gamma_{0}$; the other direction also holds: in fact, the only case in which a $\sigma$ permutes blocks of $\Gamma_{0}$ without permuting anything in $\mathcal{B}'|_{\Delta}$ is when $\Delta$ represents $k$-subsets of $\Gamma_{0}$ intersecting all parts of $\Gamma_{0}$ equally, but then there would be only one block in $\mathcal{B}'|_{\Delta}$ itself in contradiction with the fact that $|\Delta|>\frac{2}{3}|\Omega|$.
\end{proof}

This corollary divides $\Omega$ into orbits and blocks that are coarser than the original $\mathcal{B}$: this makes $n$ decrease or the block size increase, or both.

\begin{corollary}\label{co:certlargesym}
Let $|\Omega|=n$, $G\leq\mathrm{Sym}(\Omega)$ and let $\mathbf{x},\mathbf{y}:\Omega\rightarrow\Sigma$ be two strings; let $\mathcal{B}$ be a system of blocks such that $G$ acts on it as $\mathrm{Alt}(\Gamma)$ acts on $\binom{\Gamma}{k}$, where $|\Gamma|=m$ and $|\mathcal{B}|=\binom{m}{k}$. Suppose also that there exist sets $S_{\mathbf{x}},S_{\mathbf{y}}\subseteq\Gamma$ of size $>\frac{1}{2}|\Gamma|$, canonical for $\mathbf{x},\mathbf{y}$ respectively, such that for any $\sigma\in\mathrm{Alt}(S_{\mathbf{x}})$ there is an element of $\mathrm{Aut}_{G}(\mathbf{x})$ inducing $\sigma$ on $S_{\mathbf{x}}$ (and similarly for $\mathbf{y}$).

Then in time $O(n^{10})$ we can reduce the problem of determining $\mathrm{Iso}_{G}(\mathbf{x},\mathbf{y})$ to determining $4$ sets $\mathrm{Iso}_{N}(\mathbf{x},\mathbf{y}_{i})$, where $N$ induces orbits of size $\leq\frac{2}{3}|\Omega|$.
\end{corollary}

\begin{proof}
This corollary covers case \eqref{th:certlargesym} of Theorem~\ref{th:cert}.

If $\pi$ is the map going from $G$ to $\mathrm{Alt}(\Gamma)$ mentioned in the statement, define $N=\pi^{-1}(\mathrm{Alt}(\Gamma)_{(S_{\mathbf{x}})})$: we can find $N$ in time $O(n^{10})$ by Corollary~\ref{co:schreier}\eqref{co:schreierpreim}-\ref{co:schreier}\eqref{co:schreierpwstab}. Also, define $N'=\pi^{-1}(\mathrm{Alt}(\Gamma)_{S_{\mathbf{x}}})$: since $S_{\mathbf{x}}$ is canonical for $\mathbf{x}$, $\mathrm{Aut}_{G}(\mathbf{x})$ stabilizes $S_{\mathbf{x}}$ setwise, which means that it is contained inside $N'$. For any even permutation of $\Gamma$ sending $S_{\mathbf{x}}$ to $S_{\mathbf{y}}$, we can find a preimage $\tau\in G$ in time $O(n^{10})$ by Corollary~\ref{co:schreier}\eqref{co:schreierpreim}; we have
\begin{equation*}
\mathrm{Iso}_{G}(\mathbf{x},\mathbf{y})=\mathrm{Iso}_{G}(\mathbf{x},\mathbf{y}^{\tau^{-1}})\tau=\mathrm{Iso}_{N'}(\mathbf{x},\mathbf{y}^{\tau^{-1}})\tau=\mathrm{Aut}_{N'}(\mathbf{x})\mathrm{Iso}_{N}(\mathbf{x},\mathbf{y}^{\tau^{-1}})\tau,
\end{equation*}
using Lemma~\ref{le:chain}\eqref{le:chain1}, the fact that $G\tau=G$, and (by canonicity) the fact that any string isomorphism between $\mathbf{x}$ and $\mathbf{y}^{\tau^{-1}}$ must stabilize $S_{\mathbf{x}}$.

Now we have to describe $\mathrm{Aut}_{N'}(\mathbf{x})$: by the canonicity of $S_{\mathbf{x}}$, it is equal to $\mathrm{Aut}_{G}(\mathbf{x})$. Since by hypothesis $\mathrm{Alt}(S_{\mathbf{x}})$ is contained in $\mathrm{Aut}_{G}(\mathbf{x})$, there exist two elements in $\mathrm{Aut}_{G}(\mathbf{x})$ that induce two generators of $\mathrm{Alt}(S_{\mathbf{x}})$; to find them, we can take preimages $\sigma_{1},\sigma_{2}$ of these two generators in $G$ (again in time $O(n^{10})$ by Corollary~\ref{co:schreier}\eqref{co:schreierpreim}) and then determine the sets $\mathrm{Aut}_{N\sigma_{i}}(\mathbf{x})=\mathrm{Iso}_{N}(\mathbf{x},\mathbf{x}^{\sigma_{i}^{-1}})\sigma_{i}$ for $i=1,2$: any two elements $\tau_{1},\tau_{2}$ inside them will give us the whole $\mathrm{Aut}_{N'}(\mathbf{x})$, since this is $\langle A\cup\{\tau_{1},\tau_{2}\}\rangle$ for any set $A$ of generators of $\mathrm{Aut}_{N}(\mathbf{x})$. We have reduced the problem to the four problems $\mathrm{Iso}_{N}(\mathbf{x},\mathbf{y}_{i})$ with $\mathbf{y}_{1}=\mathbf{x}$, $\mathbf{y}_{2}=\mathbf{y}^{\tau^{-1}}$, $\mathbf{y}_{3}=\mathbf{x}^{\sigma_{1}^{-1}}$, $\mathbf{y}_{4}=\mathbf{x}^{\sigma_{2}^{-1}}$.

We still have to prove that $N$ has the property described in the statement. The partition $\{S_{\mathbf{x}},\Gamma\setminus S_{\mathbf{x}}\}$ can be seen as a coloured partition where $S_{\mathbf{x}}$ and $\Gamma\setminus S_{\mathbf{x}}$ are two parts of different colours (if $S_{\mathbf{x}}=\Gamma$ then the second part is empty, but this will not be a problem): examining the proof of Lemma~\ref{le:parts}, we see that each subset $\Omega_{a}$ collecting (the elements contained in blocks corresponding to) the $k$-subsets of $\Gamma$ containing $a>0$ elements of $\Gamma\setminus S_{\mathbf{x}}$ is of size $\leq\frac{2}{3}|\Omega|$; on the other hand, the blocks corresponding to $k$-subsets of $S_{\mathbf{x}}$ are stabilized by $N$ since this subgroup stabilizes $S_{\mathbf{x}}$ itself pointwise. Therefore $N$ has only orbits of size $\leq\frac{2}{3}|\Omega|$.
\end{proof}

Again, this corollary makes $n$ decrease or the block size increase (or both) by dividing $\Omega$ into orbits and blocks coarser than $\mathcal{B}$.

\begin{corollary}\label{co:certjohn}
Let $|\Omega|=n$, $G\leq\mathrm{Sym}(\Omega)$ and let $\mathbf{x},\mathbf{y}:\Omega\rightarrow\Sigma$ be two strings; let $\mathcal{B}$ be a system of blocks such that $G$ acts on it as $\mathrm{Alt}(\Gamma)$ acts on $\binom{\Gamma}{k}$, where $|\Gamma|=m$ and $|\mathcal{B}|=\binom{m}{k}$. Suppose also that, fixing the images $(y_{i})_{i=1}^{s}$ of some elements $(x_{i})_{i=1}^{s}\subseteq\Gamma$, we can find two disjoint sets $V_{1},V_{2}\subseteq\Gamma$, with $V_{2}$ divided into a system of (possibly size $1$) blocks $\mathcal{G}$ with $\binom{|\mathcal{G}|}{k'}=|V_{1}|\geq\frac{2}{3}|\Gamma|$ for some $k'\geq 2$, and a bijection between $V_{1}$ and $\binom{\mathcal{G}}{k'}$ such that each element of $G$, seen as a permutation in $\mathrm{Sym}(\mathcal{G})$, also induces the natural permutation of $V_{1}$ given by the previous identification.

Then, if $N$ is the preimage of $\mathrm{Alt}(\Gamma)_{(x_{1},\ldots,x_{s})}$ inside $G$ and $\Delta\subseteq\Omega$ is an orbit induced by $N$ of size $>\frac{2}{3}|\Omega|$, $N|_{\Delta}$ respects a system $\mathcal{B}'$ of blocks inside $\Delta$ (at least as coarse as $\mathcal{B}|_{\Delta}$), and if $M$ is the stabilizer of $\mathcal{B}'$ then $N|_{\Delta}/M\leq\mathrm{Sym}(\mathcal{G})$ (and $|\mathcal{G}|<1+\sqrt{2m}$).
\end{corollary}

\begin{proof}
This corollary covers case \eqref{th:certjohn} of Theorem~\ref{th:cert}. The focus on $N$ is due to the reduction to the problem of determining $\mathrm{Iso}_{N}(\mathbf{x},\mathbf{y}^{\sigma^{-1}})$ featured in Remark~\ref{re:multcost}, where $\sigma\in G$ is an element that sends each $x_{i}$ to $y_{i}$.

We can see $\{V_{1},V_{2},\Gamma\setminus(V_{1}\cup V_{2})\}$ as a coloured partition on $\Gamma$, where the last two parts are of size $\leq\frac{1}{3}|\Gamma|$ combined. Looking at the proof of Lemma~\ref{le:parts}, each subset $\Omega_{a}$ collecting (the elements contained in blocks corresponding to) the $k$-subsets of $\Gamma$ containing $a>0$ elements of $\Gamma\setminus V_{1}$ is of size $\leq\frac{2}{3}|\Omega|$; thus, the orbit $\Delta$ (if it exists at all) can only be one of the orbits collecting $k$-subsets of $\Gamma$ entirely contained in $V_{1}$.

An element $B\in\mathcal{B}|_{\Delta}$ corresponds to a $k$-subset $R$ of $V_{1}$ and each element of $R$ is a $k_{0}$-subset of $\mathcal{G}$; each element of $N|_{\Delta}$ induces a permutation of $\mathcal{G}$, so any two subsets $R,R'$ whose elements cover the same blocks of $\mathcal{G}$ (rather, their union does) move together under the action of $N|_{\Delta}$, i.e.\ they are in a same block of $\Delta$. A system of blocks $\mathcal{B}'$ is therefore at least as coarse as the system formed by collecting all the $B$ corresponding to the $R$ based on the same blocks of $\mathcal{G}$, which is in turn at least as coarse as $\mathcal{B}$; the image of a block $B'\in\mathcal{B}'$ is determined by the movement of the blocks of $\mathcal{G}$, since a permutation of $\mathcal{G}$ determines the new $k_{0}$-subsets of $\mathcal{G}$ represented in $V_{1}$, so $N|_{\Delta}/M\leq\mathrm{Sym}(\mathcal{G})$.

The fact that $|\mathcal{G}|<1+\sqrt{2m}$, which will be helpful in the recursion process, is evident from the hypotheses we made in the statement: since $V_{1}\subseteq\Gamma$ is in bijection with $\binom{\mathcal{G}}{k'}$ and $k'\geq 2$ we have $m\geq\binom{|\mathcal{G}|}{2}>\frac{(|\mathcal{G}|-1)^{2}}{2}$, and the inequality follows.
\end{proof}

This corollary either decreases $n$ or reduces the degree of the symmetric group that contains $G$ (as an abstract group, in the sense that we do not care about the precise action). In fact, while recursing through Cameron in this circumstance, if $G$ is not too small we will obtain a subgroup of $G$ that is $\mathrm{Alt}(\Gamma')$ for some $\Gamma'$, and $|\Gamma'|\leq 1+\sqrt{2m}$ where $m$ was the size of the old $\Gamma$.

\subsection{The algorithm, not assuming CFSG}\label{se:shortexpralgpyber}

Now we examine what the algorithm looks like when we are not assuming CFSG: the result by Cameron and Mar\'oti, which provided us with the initial crossroads to guide us in the recursion, does not hold anymore. On the other hand, the fact that the action of $G/N$ on $\mathcal{B}$ is the same as the action of $\mathrm{Alt}(\Gamma)$ on $\binom{\Gamma}{k}$ (in Theorem~\ref{th:cammar}\eqref{th:cammaralt}, Corollary~\ref{co:cammar}\eqref{co:cammaralt} and beyond) is not always essential: in many occasions the important fact is that each block of $\mathcal{B}$ corresponds to a $k$-subset of a certain $\Gamma$, but $G/N$ may act on it as some $H\leq\mathrm{Sym}(\Gamma)$, and not necessarily as $H=\mathrm{Alt}(\Gamma)$. We will see this in the next results.

We start with our new building block, a result due to Pyber \cite{Pyb93} that replaces Cameron and does not depend on CFSG.

\begin{theorem}\label{th:pyber}
Let $|\Gamma|=m$ and let $G\leq\mathrm{Sym}(\Gamma)$. Do not assume CFSG. If $G$ is primitive, then one of the following alternatives holds:
\begin{enumerate}[(a)]
\item\label{th:pybersmall} $|G|\leq m^{8\lceil 4\log_{2}m\rceil\log_{2}m}$;
\item\label{th:pybergiant} $G$ is either $\mathrm{Sym}(\Gamma)$ or $\mathrm{Alt}(\Gamma)$;
\item\label{th:pybersoj} $G$ is transitive but not doubly transitive.
\end{enumerate}
\end{theorem}

\begin{proof}
See the proof of \cite[Thm.~A]{Pyb93}.
\end{proof}

Let us tackle each of these alternatives that emerge in our determination of $\mathrm{Iso}_{G}(\mathbf{x},\mathbf{y})$. We start again with the case of $G/N$ small enough to be able to effectively use Proposition~\ref{pr:Gimprim}.

\begin{mysage}
multe_free1=25
\end{mysage}

\begin{proposition}\label{pr:pybersmall}
Let $|\Omega|=n$, $G\leq\mathrm{Sym}(\Omega)$ and let $\mathbf{x},\mathbf{y}:\Omega\rightarrow\Sigma$ be two strings; let $\mathcal{B}$ be a system of blocks preserved by $G$, and call $N$ the stabilizer of $\mathcal{B}$: suppose that there are a set $\Gamma$ of size $m$ and a bijection between $\mathcal{B}$ and $\binom{\Gamma}{k}$ (for some $k$) such that the action of $G/N$ on $\mathcal{B}$ corresponds to the action of some $H\leq\mathrm{Sym}(\Gamma)$ on $\binom{\Gamma}{k}$. Do not assume CFSG.

If $|H|\leq m^{8\lceil 4\log_{2}m\rceil\log_{2}m}$, or if $m\leq\sage{Clim_free}e^{1/\varepsilon^{2}}(\log n)^{4+\varepsilon}$, then we can reduce the problem of determining $\mathrm{Iso}_{G}(\mathbf{x},\mathbf{y})$ to determining $\leq m^{\sage{multe_free1}e^{1/\varepsilon^{2}}(\log n)^{4+\varepsilon}}$ sets of isomorphisms $\mathrm{Iso}_{N}(\mathbf{x},\mathbf{y}_{i})$, in time $O(m^{\sage{multe_free1}e^{1/\varepsilon^{2}}(\log n)^{4+\varepsilon}}n^{10})$ and at no multiplicative cost.
\end{proposition}

\begin{proof}
The proof is very similar to part of the proof of Corollary~\ref{co:cammar}, as expected: the current proposition corresponds to the route taken by Corollary~\ref{co:cammar}\eqref{co:cammarsmall}. We add that, if we know both $\Gamma$ and the bijection, it is a polynomial-time task to find out whether the conditions on $H$ are satisfied: we can calculate $|H|$ in time $O(m^{5})$ by Corollary~\ref{co:schreier}\eqref{co:schreiersize}, which will tell us if either condition is true.

\begin{mysage}
mlim_free1=5657
multe_free11=67
multe_free12=68
\end{mysage}

First, $|H|$ is always bounded by $m!\leq m^{m+\frac{1}{2}}e^{1-m}$. For $m\leq\sage{mlim_free1-1}$ we have $\left(m+\frac{1}{2}\right)\log m+1-m\leq\sage{multe_free11}\log^{3}m$, while for $m\geq\sage{mlim_free1}$ we have $4\log_{2}m>49.8$ and then $\lceil 4\log_{2}m\rceil\leq\frac{51}{50}\frac{4}{\log 2}\log m$; hence, for any $m$,
\begin{equation*}
|H|\leq m^{8\lceil 4\log_{2}m\rceil\log_{2}m} \ \ \Longrightarrow \ \ |H|\leq m^{\max\left\{\sage{multe_free11},\frac{51}{50}\frac{32}{\log^{2}2}\right\}\log^{2}m}<m^{\sage{multe_free12}\log^{2}m}.
\end{equation*}
As for $m\leq\sage{Clim_free}e^{1/\varepsilon^{2}}(\log n)^{4+\varepsilon}$, this implies easily that $|H|<m^{m}\leq m^{\sage{Clim_free}e^{1/\varepsilon^{2}}(\log n)^{4+\varepsilon}}$. Since $m\leq n$, for $\varepsilon$ small we have $\sage{multe_free12}\log^{2}m\leq\sage{Clim_free}e^{1/\varepsilon^{2}}(\log n)^{4+\varepsilon}$, so both bounds on $|H|$ can be summed up by using the unique bound $m^{\sage{multe_free1}e^{1/\varepsilon^{2}}(\log n)^{4+\varepsilon}}$. We can conclude the proof by producing all the elements of $G/N$ and working as in Proposition~\ref{pr:Gimprim}.
\end{proof}

Case \eqref{th:pybergiant} of Theorem~\ref{th:pyber} is extremely similar to the process followed in the CFSG case, as shown in the following proposition.

\begin{mysage}
dlim=16
multe_free2=ceil( alim_free+fmsoj/RIF(1/(1/10)^2)^3+loginterval(RIF(2))/RIF(1/(1/10)^2)^5 )
\end{mysage}

\begin{proposition}\label{pr:pybergiant}
Let $|\Omega|=n$, $G\leq\mathrm{Sym}(\Omega)$ and let $\mathbf{x},\mathbf{y}:\Omega\rightarrow\Sigma$ be two strings; let $\mathcal{B}$ be a system of blocks preserved by $G$, and call $N$ the stabilizer of $\mathcal{B}$: suppose that there are a set $\Gamma$ of size $m$ and a bijection between $\mathcal{B}$ and $\binom{\Gamma}{k}$ (for some $k$) such that the action of $G/N$ on $\mathcal{B}$ corresponds to the action of $H=\mathrm{Sym}(\Gamma),\mathrm{Alt}(\Gamma)$ on $\binom{\Gamma}{k}$. Do not assume CFSG.

If $m>\sage{Clim_free}e^{1/\varepsilon^{2}}(\log n)^{4+\varepsilon}$, then we reduce the problem of determining $\mathrm{Iso}_{G}(\mathbf{x},\mathbf{y})$ to one of the following:
\begin{enumerate}[(a)]
\item\label{pr:pybergiantlargesym} determining $\leq 8$ sets $\mathrm{Iso}_{N'}(\mathbf{x},\mathbf{y}_{i})$, where $N'$ divides $\Omega$ into orbits of size $\leq\frac{2}{3}|\Omega|$;
\item\label{pr:pybergiantsoj} determining $\leq m^{\sage{multe_free2}e^{1/\varepsilon^{2}}(\log n)^{4+\varepsilon}}$ sets $\mathrm{Iso}_{N'}(\mathbf{x},\mathbf{y}_{i})$, where $N'$ divides $\Omega$ into a system of orbits and/or blocks $\mathcal{B}'$ (at least as coarse as $\mathcal{B}$) such that if there is an orbit $\Delta$ of size $>\frac{2}{3}|\Omega|$ then either
\begin{enumerate}[(\ref{pr:pybergiantsoj}1)]
\item\label{pr:pybergiantpart} $\mathcal{B}'|_{\Delta}$ is nontrivial and strictly coarser than $\mathcal{B}|_{\Delta}$, with stabilizer of $\mathcal{B}'|_{\Delta}$ equal to the block stabilizer of the large colour of $\Gamma$ (in the sense of Corollary~\ref{co:certpart}), or
\item\label{pr:pybergiantjohn} if $M$ is the stabilizer of $\mathcal{B}'|_{\Delta}$, $N'|_{\Delta}/M$ acts on $\mathcal{B}'|_{\Delta}$ as some $H'\leq\mathrm{Sym}(\Gamma')$ acts on $\binom{\Gamma'}{k'}$ with $|\Gamma'|<1+\sqrt{2m}$.
\end{enumerate}
\end{enumerate}
The time necessary for this reduction is the cost of $\frac{1}{2}m^{2a}naa!$ calls of the whole algorithm for strings of length $\leq\frac{n}{a}$ where  $a\in(\sage{alim_freelow},\sage{alim_free})\cdot e^{1/\varepsilon^{2}}(\log n)^{4+\varepsilon}$, plus some additional time $O(m^{3a}n^{11})$.
\end{proposition}

\begin{proof}
First, in the case of $H=\mathrm{Sym}(\Gamma)$ we can reduce the problem to $2$ sets with $H=\mathrm{Alt}(\Gamma)$. Now we are exactly in the case described in Corollary~\ref{co:cammar}\eqref{co:cammaralt}. We can retrace all the steps from Theorem~\ref{th:cert} to Corollary~\ref{co:certjohn}, this time using the CFSG-free versions of the results in \S\ref{se:majorou}, and the results correspond to one of the final situations thereby reached: case \eqref{pr:pybergiantlargesym} corresponds to Corollary~\ref{co:certlargesym} (where $4$ becomes $8$ because of the aforementioned reduction from $\mathrm{Sym}$ to $\mathrm{Alt}$), case \eqref{pr:pybergiantpart} corresponds to Corollary~\ref{co:certpart}, and case \eqref{pr:pybergiantjohn} corresponds to Corollary~\ref{co:certjohn}.

We need only to justify how to obtain the action in part \eqref{pr:pybergiantjohn} rather than only a bound on the degree of $N|_{\Delta}/M$ like in Corollary~\ref{co:certjohn} (as we observed, this stronger statement is necessary for the recursion, given the unavailability of Cameron).

Let us start with the first problem. Following the reasoning up to Corollary~\ref{co:certjohn}, we ended up finding two disjoint sets $V_{1},V_{2}\subseteq\Gamma$ and a partition $\mathcal{G}$ of $V_{2}$ that respect the various hypotheses mentioned in the corollary, and in its proof we find a system of blocks $\mathcal{B}'|_{\Delta}$ on an orbit $\Delta$ of size $>\frac{2}{3}|\Omega|$ (if such an orbit exists) such that the action of $N|_{\Delta}$ is induced by the permutations of $\mathcal{G}$, up to the stabilizer of the system. If $k=1$, $\mathcal{B}$ corresponds to $\Gamma$ itself: therefore $\Delta$ of size $>\frac{2}{3}|\Omega|$ must correspond to $V_{1}$ itself, and by hypothesis the permutations of $\mathcal{G}$ induce permutations of $V_{1}$ in a way that respects the bijection $V_{1}\leftrightarrow\binom{\mathcal{G}}{k'}$ ($\mathcal{G}$ is then the sought $\Gamma'$). If $k\geq 2$, we can use Lemma~\ref{le:johnjohn} to prove that $\Delta$ is further split into blocks that are strictly coarser than $\mathcal{B}$: in that lemma, we use $\Gamma',\Gamma,\mathcal{B}$ to refer in this situation to $\mathcal{G},V_{1},\mathcal{B}|_{\Delta}$ respectively; we only have to show that the bounds on $|\mathcal{G}|$ hold. If $m>\sage{Clim}\log^{2}n$, by Remark~\ref{re:thism} we have $m\geq\sage{mlim}$; $|V_{1}|\geq\frac{2}{3}m$, so that $|V_{1}|\geq\sage{ceil(2/3*mlim)}$: whatever will be our choice of $k'$, we have $\sage{ceil(2/3*mlim)}\leq\binom{|\mathcal{G}|}{k'}\leq\binom{|\mathcal{G}|}{\lfloor\frac{1}{2}|\mathcal{G}|\rfloor}$, hence $|\mathcal{G}|\geq\sage{mplim}$.

Finally, let us obtain the exponent in part \eqref{pr:pybergiantsoj} and the value of $a$. The interval of $a$ is taken directly from Proposition~\ref{pr:computcert}\eqref{pr:computcertfree}. As for the exponent, we notice that exactly as in Theorem~\ref{th:cert} we still have $d\leq a$ (with $d$ as in Lemma~\ref{le:dtrans}\eqref{le:dtransfree}), so that the multiplicative cost is still $m^{a+\sage{fmsoj}\log m}$. For our choice of $a$, our bounds $\sage{Clim_free}e^{1/\varepsilon^{2}}(\log n)^{4+\varepsilon}<m\leq n$, and $\varepsilon$ small enough, we can bound this cost as in the statement (remember that we also have a possible multiplication by $2$, from the reduction in the case of $H=\mathrm{Sym}(\Gamma)$). The additive cost is the same as in Theorem~\ref{th:cert}.
\end{proof}

Finally, we treat case \eqref{th:pybersoj} of Theorem~\ref{th:pyber}, whose procedure is a somewhat shortened version of the one covered in the previous proposition.

\begin{proposition}\label{pr:pybersoj}
Let $|\Omega|=n$, $G\leq\mathrm{Sym}(\Omega)$ and let $\mathbf{x},\mathbf{y}:\Omega\rightarrow\Sigma$ be two strings; let $\mathcal{B}$ be a system of blocks preserved by $G$, and call $N$ the stabilizer of $\mathcal{B}$: suppose that there are a set $\Gamma$ of size $m$ and a bijection between $\mathcal{B}$ and $\binom{\Gamma}{k}$ (for some $k$) such that the action of $G/N$ on $\mathcal{B}$ corresponds to the action of some $H\leq\mathrm{Sym}(\Gamma)$ on $\binom{\Gamma}{k}$. Do not assume CFSG.

If $m>\sage{Clim_free}e^{1/\varepsilon^{2}}(\log n)^{4+\varepsilon}$ and $H$ is transitive but not doubly transitive, then in time $O(m^{14})$ we reduce the problem of determining $\mathrm{Iso}_{G}(\mathbf{x},\mathbf{y})$ to determining $\leq m^{\sage{fmsoj}\log m}$ sets $\mathrm{Iso}_{N'}(\mathbf{x},\mathbf{y}_{i})$ where $N'$ divides $\Omega$ into a system of orbits and/or blocks $\mathcal{B}'$ (at least as coarse as $\mathcal{B}$) such that if there is an orbit $\Delta$ of size $>\frac{2}{3}|\Omega|$ then either
\begin{enumerate}[(a)]
\item\label{pr:pybersojpart} $\mathcal{B}'|_{\Delta}$ is nontrivial and strictly coarser than $\mathcal{B}|_{\Delta}$, with stabilizer of $\mathcal{B}'|_{\Delta}$ equal to the block stabilizer of the large colour of $\Gamma$ (in the sense of Corollary~\ref{co:certpart}), or
\item\label{pr:pybersojjohn} if $M$ is the stabilizer of $\mathcal{B}'|_{\Delta}$, $N'|_{\Delta}/M$ acts on $\mathcal{B}'|_{\Delta}$ as some $H'\leq\mathrm{Sym}(\Gamma')$ acts on $\binom{\Gamma'}{k'}$ with $|\Gamma'|<1+\sqrt{2m}$.
\end{enumerate}
\end{proposition}

\begin{proof}
If $H$ is transitive but not doubly transitive, we can determine the nontrivial orbits of the action of $H$ on $\binom{\Gamma}{2}$ in time $O(m^{6})$ by Lemma~\ref{le:orbblo}; giving to each orbit its own colour, we can make $\binom{\Gamma}{2}$ into a coherent configuration in time $O(m^{10}\log m)$ (mostly due to Weisfeiler-Leman, see \cite[\S\S 2.3-2.5]{Hel19}): the result would be a nontrivial homogeneous coherent configuration, where homogeneity is consequence of the fact that this is a canonical process and $H$ moves every point of $\Gamma$ to any other, so that we are unable to distinguish them with different colours.

Now we can use SoJ directly. We use Proposition~\ref{pr:computsoj}, where from the costs we can remove the exponent $b$ (since we do not perform the Design Lemma).

The shape of the action of $N'|_{\Delta}/M$ on $\mathcal{B}'|_{\Delta}$ in part \eqref{pr:pybersojjohn} is again proved as in part \eqref{pr:pybergiantjohn} of Proposition~\ref{pr:pybergiant}, i.e.\ resorting to Lemma~\ref{le:johnjohn}.
\end{proof}

All these cases reduce to some sort of recursion with lower parameters, either by decreasing $n$ or $m$ or increasing the block size. This works exactly as in the CFSG case.

\section{Main theorem: proof}\label{se:shortexprcost}

We are at last ready to prove Theorem~\ref{th:shortexpr}.

The group-theoretic results to which we keep returning in our recursions are Theorem~\ref{th:cammar} in the CFSG case and Theorem~\ref{th:pyber} in the CFSG-free case; we have already declared this multiple times, but we repeat it here (now with references, though): except for exiting through the base cases given in Remark~\ref{re:nsmall} and Proposition~\ref{pr:base} and for breaking down $\Omega$ into smaller orbits through Proposition~\ref{pr:Gintrans}, the only other alternatives are that on a large chunk of $\Omega$ either the system of blocks $\mathcal{B}$ on which we are working becomes coarser and coarser (the conclusion featured in Corollary~\ref{co:certpart}, Proposition~\ref{pr:pybergiant}\eqref{pr:pybergiantpart} and Proposition~\ref{pr:pybersoj}\eqref{pr:pybersojpart}) or the group in which we are operating is contained in a symmetric group of degree smaller and smaller (the conclusion featured in Corollary~\ref{co:certjohn}, Proposition~\ref{pr:pybergiant}\eqref{pr:pybergiantjohn} and Proposition~\ref{pr:pybersoj}\eqref{pr:pybersojjohn}).

\begin{proof}[Proof of Thm.~\ref{th:shortexpr}]
There are several tasks to accomplish: we need to analyze the possible passages mentioned above and see that they fit the description given in terms of ($\mathcal{C}$1)-($\mathcal{C}$2)-($\mathcal{C}$3), and that the final base cases fit ($\mathcal{A}$), and we need to estimate their contribution in terms of both the multiplicative cost (which will lead us to a bound on the number of atomic elements) and additive cost (which will yield the total runtime).

To determine the multiplicative cost of the procedure, we start in medias res. We are working on a certain orbit $\Delta$ of $\Omega$, of size $|\Delta|=n'\leq n$, divided into a system of blocks $\mathcal{B}$, of size $|\mathcal{B}|=r\leq n'$, such that the group $G/N$ permuting the blocks is isomorphic to a subgroup of $\mathrm{Sym}(m)$, of degree $m\leq r$. We call $M(n',r,m)$ (an upper bound on) the multiplicative cost that we incur from this moment until we manage to make each block into an orbit of its own. Call $T(n',r,m)$ the intermediate time cost, in an analogous fashion as we did with $M(n',r,m)$; we also suppose that $T(n',r,m)$ includes the cost of performing Proposition~\ref{pr:Gintrans} on the resulting orbits, so as to cover the time spent to bridge one intermediate problem to the next one.

The proof is articulated in the following main steps. 
\begin{enumerate}[(1)]
\item From the already known passages we delineate a handful of ``actions'' and the reduction they entail on $M(n',r,m)$; note that here we are using the word ``action'' \textit{not} in a mathematical sense, but in the everyday meaning of ``something done purposefully to accomplish a certain end''. This step gives us a series of conditions that our function $M$ must respect in order to work.
\item We choose $M$ and show that it is compatible with the previous conditions coming from the actions; then $M(n,n,n)$ by definition turns out to be a bound on the multiplicative cost incurred throughout the whole algorithm.
\item We translate actions into ($\mathcal{C}$1)-($\mathcal{C}$2)-($\mathcal{C}$3) and end-cases into ($\mathcal{A}$), and use $M(n,n,n)$ to bound the number of atomic elements.
\item We refine the computations of the second part to tackle $T(n',r,m)$.
\end{enumerate}
For the sake of notation, we are going to perform our computations by bounding $\log M$ instead of $M$, so that the focus will be on the exponents of the quantities involved.

\smallskip
(1) \textit{Description of the actions.}

\begin{mysage}
K1_cfsg=Clim
K1_free=multe_free1
K2_cfsg=multe_cfsg
K2_free=multe_free2
\end{mysage}

The first action that is possible to perform, following from Corollary~\ref{co:cammar}\eqref{co:cammarsmall} and Proposition~\ref{pr:pybersmall}, is to directly pass to the stabilizer of the system, thus making each block into an orbit: this concludes the calculation of $M$ with no reduction, and it costs at most $\sage{Clim}\log m\log^{2}n'$ in the CFSG case and $\sage{multe_free1}e^{1/\varepsilon^{2}}\log m(\log n')^{4+\varepsilon}$ in the CFSG-free case; these are direct lower bounds for $\log M(n',r,m)$, therefore
\begin{equation}\label{eq:act1}
\log M(n',r,m)\geq K_{1}\log m(\log n')^{e_{1}}
\end{equation}
for $(K_{1},e_{1})=(\sage{K1_cfsg},2),(\sage{K1_free}e^{1/\varepsilon^{2}},4+\varepsilon)$ appropriately.

For notational simplicity, let us set $X=\sage{mlim}$ for the CFSG case and $X=\sage{mlim_free}e^{1/\varepsilon^{2}}$ for the CFSG-free case: these are the values we have already encountered many times, and they separate small and large values of $m,n$ (see Remark~\ref{re:thism} in particular). If either $n'$ or $m$ is smaller than $X$ we are using the first action, so for the other actions we can assume otherwise.

The second action, following from Corollary~\ref{co:certlargesym} and Proposition~\ref{pr:pybergiant}\eqref{pr:pybergiantlargesym} and (in case there are only orbits of size $\leq\frac{2}{3}|\Omega|$) from Corollaries~\ref{co:certpart}-\ref{co:certjohn} and Propositions~\ref{pr:pybergiant}\eqref{pr:pybergiantsoj}-\ref{pr:pybersoj}, consists in reducing $n'$ (and consequently $r$) by a fraction at least as small as $\frac{2}{3}$. This costs at most $K_{2}\log m(\log n')^{e_{2}}$, where $(K_{2},e_{2})=(\sage{K2_cfsg},1)$ assuming CFSG and $(K_{2},e_{2})=(\sage{K2_free}e^{1/\varepsilon^{2}},4+\varepsilon)$ without CFSG: for our bounds on $m,n'$ (and for $\varepsilon$ small), these are the largest expenses, coming from Theorem~\ref{th:cert}\eqref{th:certsoj} and Propositions~\ref{pr:pybergiant}\eqref{pr:pybergiantsoj} respectively. Hence
\begin{equation}\label{eq:act2}
\log M(n',r,m)\geq K_{2}\log m(\log n')^{e_{2}}+\log M\left(\frac{2}{3}n',\frac{2}{3}r,m\right).
\end{equation}

The third action, following (in case there is an orbit of size $>\frac{2}{3}|\Omega|$) from Corollary~\ref{co:certpart} and Propositions~\ref{pr:pybergiant}\eqref{pr:pybergiantpart}-\ref{pr:pybersoj}\eqref{pr:pybersojpart}, creates a new system of blocks strictly coarser than the original $\mathcal{B}$, at a cost of at most $K_{2}\log m(\log n')^{e_{2}}$: $(K_{2},e_{2})$ is as in the previous action, as the largest expenses originate in the same results. What happens is, we have first to work on the coarser system, then after we have stabilized each coarser block we have to work on each one of them as the new orbit and the finer blocks as the new system; since the stabilizer of coarser blocks coincides with some block stabilizer of $\Gamma$, we also get $m',\frac{m}{m'}$ instead of $m$ in the two steps, for some $2\leq m'\leq\frac{m}{2}$. The bound on $\log M(n',r,m)$ given by this action is
\begin{equation}\label{eq:act3}
\log M(n',r,m) \geq K_{2}\log m(\log n')^{e_{2}}+\log M(n',r',m')+\log M\left(\frac{n'}{r'},\frac{r}{r'},\frac{m}{m'}\right),
\end{equation}
where $2\leq r'\leq\frac{r}{2}$ is the size of the coarser system.

The fourth action, following (in case there is an orbit of size $>\frac{2}{3}|\Omega|$) from Corollary~\ref{co:certjohn} and Propositions~\ref{pr:pybergiant}\eqref{pr:pybergiantjohn}-\ref{pr:pybersoj}\eqref{pr:pybersojjohn}, reduces the degree of the minimal symmetric group containing $G$, at a cost of at most $K_{2}\log m(\log n')^{e_{2}}$ ($(K_{2},e_{2})$ as in the second and third actions); therefore,
\begin{equation}\label{eq:act4}
\log M(n',r,m)\geq K_{2}\log m(\log n')^{e_{2}}+\log M(n',r,1+\sqrt{2m}).
\end{equation}

\smallskip
(2) \textit{Choice of function $M$.}

Now let us prove that
\begin{equation}\label{eq:intermult}
\log M(n',r,m)=(\log n')^{e_{2}+1}(a\log m+b\log r)
\end{equation}
satisfies the four conditions for some appropriate constants $a,b$.

Since $m\leq r$ and $e_{1}\leq e_{2}+1$, in order to have \eqref{eq:act1} we have simply to ask $a+b\geq K_{1}$. Recall that for the other actions we can assume $m,n'\geq X$.

For $n'\geq X$ and $e_{2}\geq 1$ we have $\left(\log\left(\frac{2}{3}n'\right)\right)^{e_{2}+1}<(\log n')^{e_{2}+1}-\frac{3}{4}(\log n')^{e_{2}}$ (for both values of $X$), so
\begin{align*}
 & \ K_{2}\log m(\log n')^{e_{2}}+\left(\log\left(\frac{2}{3}n'\right)\right)^{e_{2}+1}\left(a\log m+b\log\left(\frac{2}{3}r\right)\right)\\
 < & \ (\log n')^{e_{2}+1}(a\log m+b\log r)+(\log n')^{e_{2}}\left(K_{2}\log m-\frac{3}{4}(a\log m+b\log r)\right),
\end{align*}
and since $m\leq r$ in order to have \eqref{eq:act2} it is sufficient to ask $\frac{3}{4}(a+b)>K_{2}$. For \eqref{eq:act3}, using $\left(\log\frac{n'}{r'}\right)^{e_{2}+1}<(\log n')^{e_{2}+1}-\log r'(\log n')^{e_{2}}$ and $\log\frac{m}{m'}\geq\log 2$ the sufficiency of \eqref{eq:intermult} in this case is implied by
\begin{equation}\label{eq:act3suff}
f(\log r')=b\log^{2}r'-\left(a\log 2+b\log r\right)\log r'+K_{2}\log m\leq 0.
\end{equation}
The function $f(x)$ in the interval $[\log 2,\log r-\log 2]$ has its maximum in $x=\log 2$, being a quadratic polynomial with the minimum in $x=\frac{1}{2}\log r+\frac{a\log 2}{2b}>\frac{1}{2}\log r$; evaluating $f(\log 2)$ and recalling that $X\leq m\leq r$, \eqref{eq:act3suff} is in turn consequence of
\begin{equation}\label{eq:act3b}
b\geq\frac{K_{2}\log m}{\log 2(\log r-\log 2)}-\frac{K_{2}\log 2}{\log r-\log 2}a \ \ \Longleftarrow \ \ b\geq\frac{K_{2}\log X}{\log 2\log(X/2)}.
\end{equation}

\begin{mysage}
expact4true=loginterval(1+sqrtinterval(2*RIF(mlim)))/loginterval(RIF(mlim))
expact4=roundup(expact4true,5)
invact4=round(1-expact4,ndigits=5)
act4=roundup(1/(invact4*loginterval(RIF(mlim))),5)
\end{mysage}

To have \eqref{eq:act4}, we notice that $1+\sqrt{2m}<m^{\sage{expact4}}$ for $m\geq X$ (for both values of $X$); then,
\begin{align*}
& \ (\log n')^{e_{2}+1}(a\log m+b\log r) \\
\geq & \ K_{2}\log m(\log n')^{e_{2}}+(\log n')^{e_{2}+1}(\sage{expact4}a\log m+b\log r)
\end{align*}
means $a\geq\frac{K_{2}}{\sage{invact4}\log n'}$, so that $a\geq\sage{act4}K_{2}\geq\frac{K_{2}}{\sage{invact4}\log X}$ is enough to satisfy \eqref{eq:act4}.

\begin{mysage}
a_cfsgtrue=act4*K2_cfsg
a_cfsg=roundup(a_cfsgtrue,5)
b_cfsgtrue=K2_cfsg*loginterval(RIF(mlim))/(loginterval(RIF(2))*loginterval(RIF(mlim/2)))
b_cfsg=roundup(b_cfsgtrue,5)
main_cfsgtrue=round(a_cfsg+b_cfsg,ndigits=5)
ab_free=K1_free
main_freetrue=ab_free
\end{mysage}

Putting together these conditions and considering our $K_{1},K_{2}$, it turns out that $a=\sage{a_cfsg}$ and $b=\sage{b_cfsg}$ with CFSG and $a=b=\frac{\sage{ab_free}}{2}e^{1/\varepsilon^{2}}$ without CFSG are suitable choices for \eqref{eq:intermult}. The multiplicative cost of the whole algorithm is bounded by $M(n,n,n)$; thus we conclude that the multiplicative cost is bounded by
\begin{align}\label{eq:multfinal}
n^{\sage{main_cfsgtrue}\log^{2}n} & \text{ with CFSG,} & n^{\sage{main_freetrue}e^{1/\varepsilon^{2}}(\log n)^{5+\varepsilon}} & \text{ without CFSG.}
\end{align}

\smallskip
(3) \textit{Reduction to $\mathrm{(\mathcal{A})}$-$\mathrm{(\mathcal{C}1)}$-$\mathrm{(\mathcal{C}2)}$-$\mathrm{(\mathcal{C}3)}$.}

Now that we have bounded the multiplicative cost, let us focus now on the actions themselves, in order to be able to describe the various stages as one among ($\mathcal{A}$)-($\mathcal{C}$1)-($\mathcal{C}$2)-($\mathcal{C}$3) and to use $M(n,n,n)$ for the computation of the number of atomic elements.

The first action entails firstly a reduction of the problem of determining the set $\mathrm{Iso}_{G}(\mathbf{x},\mathbf{y})$ to a collection of $\mathrm{Iso}_{N}(\mathbf{x},\mathbf{y}_{i}^{\sigma_{i}^{-1}})\sigma_{i}$ whose union is the original set, as seen in Proposition~\ref{pr:Gimprim} or Remark~\ref{re:multcost}: the way this union is performed corresponds precisely to ($\mathcal{C}$1), and the number of subproblems is equal to the multiplicative cost incurred during this action; then, each stabilized block becomes an orbit of its own, in a reduction that corresponds to the situation described in ($\mathcal{C}$2) (see Proposition~\ref{pr:Gintrans}). This passage does not feature any multiplicative cost, but it does multiply the number of atomic elements at the end: however, since we have simply $r$ blocks, the contribution of ($\mathcal{C}$2) here, and indeed the contribution of any nested series of ($\mathcal{C}$2) acting throughout the entire process of solving the intermediate problem with parameters $(n',r,m)$, is at most $r$.

The second action features a reduction of $\Omega$ to orbits of size at most $\frac{2}{3}|\Omega|$; this can happen in two different ways. In the case of Corollaries~\ref{co:certpart}-\ref{co:certjohn} and Propositions~\ref{pr:pybergiant}\eqref{pr:pybergiantsoj}-\ref{pr:pybersoj}, after having fixed the image of a certain number of points at a multiplicative cost we find orbits of such size, and then we examine each orbit singularly: this is exactly as in the previous case, where each passage consists in using ($\mathcal{C}$1) and ($\mathcal{C}$2), and the bounds on the atomic element multiplication are as above. In the case of Corollary~\ref{co:certlargesym} and Proposition~\ref{pr:pybergiant}\eqref{pr:pybergiantlargesym}, we are in a situation where
\begin{equation*}
\mathrm{Iso}_{G}(\mathbf{x},\mathbf{y})=\langle\mathrm{Aut}_{N}(\mathbf{x}),\tau_{1},\tau_{2}\rangle\tau'\tau,
\end{equation*}
where $\tau'\in\mathrm{Iso}_{N}(\mathbf{x},\mathbf{y}^{\tau^{-1}})$ (to use the notation of the corollary); this corresponds to ($\mathcal{C}$3), and despite the multiplication cost being at most $4$ or $8$, there is no actual growth in the number of atomic elements through this case.

The third and the fourth action create respectively (on the large orbit) a strictly coarser system of blocks and a bijection on a permutation subgroup of strictly smaller degree: this happens at a certain multiplicative cost, that corresponds to a passage of the form shown in ($\mathcal{C}$1) and multiplies the atomic elements by the same quantity.

The various actions, as we already said, decrease at least one of the three parameters $n,r,m$, and when $r,m$ become too small $n$ itself diminishes through the use of the first action: hence, the procedure eventually stops when $n=1$, the trivial case of Remark~\ref{re:nsmall}. There is also a second way to stop the algorithm, and that is Proposition~\ref{pr:base}: both cases correspond to the atom ($\mathcal{A}$). The reduction to ($\mathcal{A}$)-($\mathcal{C}$1)-($\mathcal{C}$2)-($\mathcal{C}$3) has been proved; the actual writing of the expression is done following the proofs of Proposition~\ref{pr:Gimprim} (for ($\mathcal{C}$1)), Proposition~\ref{pr:Gintrans} (for ($\mathcal{C}$2)) and Corollary~\ref{co:certlargesym} (for ($\mathcal{C}$3)). The number of atomic elements, by the reasoning above, is bounded by
\begin{align*}
n\cdot n^{\sage{main_cfsgtrue}\log^{2}n} & <n^{\sage{main_cfsg}\log^{2}n} & & \text{with CFSG}, \\
n\cdot n^{\sage{main_freetrue}e^{1/\varepsilon^{2}}(\log n)^{5+\varepsilon}} & <n^{\sage{main_free}e^{1/\varepsilon^{2}}(\log n)^{5+\varepsilon}} & & \text{without CFSG},
\end{align*}
since its intermediate multiplication is bounded by $rM(n',r,m)$, and we are done.

\smallskip
(4) \textit{Runtime.}

Finally, let us tackle the runtime; we start at the end, this time. We have already proved that there are at most $n^{K\log^{e}n}$ atomic elements constituting the expression, and by Remark~\ref{re:nsmall} and Proposition~\ref{pr:base} we can treat each one in time $O(n^{6})$, so the bound on the runtime covers this final stage; now we go back to the analysis of the recursion process that leads to it.

\begin{mysage}
lowlognu=rounddown(loginterval(loginterval(RIF(mlim))),5)
leewayexp=floor(5*RIF(lowlognu)*loginterval(RIF(mlim))^2/loginterval(RIF(10)))
\end{mysage}

Call $T(n',r,m)$ the intermediate time cost, in an analogous fashion as we did with $M(n',r,m)$; most of the computations for $M$ also hold for $T$, but we have to verify that the \textit{added} time does not disrupt the final constants coming from our multiplicative reasoning: we also suppose that $T(n',r,m)$ includes the cost of performing Proposition~\ref{pr:Gintrans} on the resulting orbits, so as to cover the time spent to bridge one intermediate problem to the next one. For the first action, the bound is as in Corollary~\ref{co:cammar}\eqref{co:cammarsmall} and Proposition~\ref{pr:pybersmall}, with the addition of the cost for the reduction to single orbits:
\begin{equation*}
T(n',r,m)=O(m^{K_{1}(\log n')^{e_{1}}}n'^{10}+n'^{11}).
\end{equation*}
As for the other three actions, let us start by working on the additive cost first; recall that henceforth $n'\geq r\geq m\geq X$. The highest additive cost is featured in Theorem~\ref{th:cert} and Proposition~\ref{pr:pybergiant} and it involves the use of the runtime itself (for smaller $n'$); supposing that we want to show that it is sufficient to ask $T(n',r,m)=O(e^{(\log n')^{e_{2}+1}(a\log m+b\log r)}n'^{11})$, this cost is of order
\begin{equation}\label{eq:addt}
\frac{1}{2}m^{2\nu}n'\nu\nu!\cdot e^{(\log\frac{n'}{\nu})^{e_{2}+1}(a\log m+b\log r)}\frac{n'^{11}}{\nu^{11}}+2m^{3\nu}n'^{11},
\end{equation}
where $\nu=\alpha(\log n')^{e_{2}}$ for some $\alpha\in\left(\sage{alim_low},\sage{alim}\right)$ with CFSG and $\alpha\in(\sage{alim_freelow},\sage{alim_free})\cdot e^{1/\varepsilon^{2}}$ without CFSG. Notice that we write $2m^{3\nu}n'^{11}$ (i.e.\ with a $2$ in front) in order to absorb the successive smaller costs, such as the $n'^{11}$ from Proposition~\ref{pr:Gintrans}, the $n'^{10}$ from Corollary~\ref{co:certlargesym} and the $m^{14}$ from Proposition~\ref{pr:pybersoj}. For $a,b\geq 5$, it is easy to prove that the first addend of \eqref{eq:addt} is larger than the second: say for example $n'>4$, $\nu\nu!>1$ and $e^{(\log\frac{n'}{\nu})^{e_{2}+1}(a\log m+b\log r)}>e^{\frac{1}{3}\log^{2}n'(a\log m+b)}=m^{\frac{a}{3}\log^{2}n'}n'^{\frac{b}{3}\log n'}>m^{\nu}(2\nu)^{11}$. Now let us bound the first addend (without $\frac{1}{2}$); its logarithm is
\begin{align*}
 & \ 2\nu\log m+\log(n'\nu\nu!)+\left(\log\frac{n'}{\nu}\right)^{e_{2}+1}(a\log m+b\log r)+\log\frac{n'^{11}}{\nu^{11}} \\
 < & \ 2\alpha(\log n')^{e_{2}}\log m+\log n'+\log m+\alpha(\log n')^{e_{2}}\log m \\
 & \ +(\log n')^{e_{2}+1}(a\log m+b\log r)-\sage{lowlognu}(\log n')^{e_{2}}(a\log m+b\log r)+\log n'^{11} \\
 < & \ (\log n')^{e_{2}+1}(a\log m+b\log r)+\log n'^{11}-\sage{lowlognu}b(\log n')^{e_{2}}\log r,
\end{align*}
using $\left(\log\frac{n'}{\nu}\right)^{e_{2}+1}<(\log n')^{e_{2}+1}-(\log n')^{e_{2}}\log\nu$ for $e_{2}\geq 1$ and $\sage{lowlognu}<\log\nu<\log m$, and noting that the negative $(\log n')^{e_{2}}\log m$ term absorbs the smaller $\log n',\log m,(\log n')^{e_{2}}\log m$ positive terms for $3\alpha+2<2a$. Therefore for example $b\geq 5$ gives us already enough leeway:
\begin{equation*}
e^{-\sage{lowlognu}b(\log n')^{e_{2}}\log r}<10^{-\sage{leewayexp}}.
\end{equation*}

Now that the additive cost is accounted for, we continue with the multiplicative one. Since we want to prove that a quantity multiplied by $n'^{11}$ is larger than its partial version multiplied by some fraction of $n'^{11}$, we can just ignore this polynomial cost. For the second action, we exploit the already existing margin left out before: $\left(\log\left(\frac{2}{3}n'\right)\right)^{e_{2}+1}<(\log n')^{e_{2}+1}-\left(\frac{3}{4}+\frac{3}{100}\right)(\log n')^{e_{2}}$, and for $a+b\geq 1$ we are left with a constant of
\begin{equation*}
e^{-\frac{3}{100}(\log n')^{e_{2}}(a\log m+b\log r)}<\frac{1}{4}
\end{equation*}
in front of this part of the runtime. For the third action, if $b$ is as on the right side of \eqref{eq:act3b}, we can use $\left(1+\frac{1}{100000}\right)b$ as the new coefficient and going through \eqref{eq:act3suff} we can cut ourselves a margin of
\begin{equation*}
e^{-\frac{b}{100000}(\log n')^{e_{2}}\log r'\log\frac{r}{r'}}\leq e^{-\frac{K_{2}}{100000}(\log n')^{e_{2}}\log X}<\frac{49}{50}.
\end{equation*}
The fourth action is treated in the same way: putting $\left(1+\frac{1}{100000}\right)a$ we carve out a $\frac{49}{50}$ constant as well. This shows that we can take the same coefficient $a,b$ as before multiplied by $1+\frac{1}{100000}$, because $\frac{49}{50}+10^{-\sage{leewayexp}}<1$; also, thanks to
\begin{align*}
n^{\sage{main_cfsgtrue}\left(1+\frac{1}{100000}\right)\log^{2}n} & <n^{\sage{main_cfsg}\log^{2}n}, \\
n^{\sage{main_freetrue}\left(1+\frac{1}{100000}\right)e^{1/\varepsilon^{2}}(\log n)^{5+\varepsilon}} & <n^{\sage{main_free}e^{1/\varepsilon^{2}}(\log n)^{5+\varepsilon}},
\end{align*}
we achieve the bounds we wanted in the two cases for the runtime, too.

\smallskip
The theorem is proved.
\end{proof}

\section{Concluding remarks}

It must be noted that the difference between the exponents for the CFSG and the CFSG-free case in not a consequence of the different use of group-theoretic results to produce a suitable recursion (Theorems~\ref{th:cammar} and \ref{th:pyber} respectively): they make the algorithm different in the two cases, that is true, but the different expense lies elsewhere. What is important in this respect is the theoretic tool that allows the recursion in Theorem~\ref{th:cert} and Proposition~\ref{pr:pybergiant}, and that gives for us a different number of calls to the algorithm for shorter strings. In the local certificates procedure in Babai's algorithm, one important detail is that a certain epimorphism $G\rightarrow\mathrm{Alt}(k)$ for $G\leq\mathrm{Sym}(n)$ primitive is guaranteed to be an isomorphism, and this is ensured for $k=\Omega(\log n)$ with a proof relying on CFSG (see \cite[Lemma 8.3.1]{Bab16a} \cite[Lemme 4.1]{Hel19}), but only for $k=\Omega(\log n)^{4+\varepsilon}$ without CFSG (see \cite[Lemma 12]{Pyb16}, where $\Omega(\log^{5}n)$ is used). Consequently the algorithm is still performing the same subroutines, but the tuples on which we want to build the certificates need to be larger, leading to the loss of efficiency that we witness.

The constants are likely improvable, if one were to analyze with greater care the routines. We have been quite accurate, but we have not really aimed at obtaining the best possible constant, especially in the CFSG-free case: as our position is to consider CFSG as a theorem, the analysis of the CFSG-free procedure is more of a question of method, especially given the way the main theorem is applied in \cite[\S 6]{Don20}.

\begin{center}
***
\end{center}

In truth, the origin of the whole analysis performed in here lay originally in trying to find whether we could easily arrive to an improvement of Babai's algorithm that would gets us to a $n^{O(\log n)}$ runtime, or, if not, to point out where exactly the bottleneck was and why.

It is clear, to the attentive reader of these pages, that the obstacle does not lie in the ``interstitial reasoning'' as we called it at the start. We have performed our analysis burdened with multiplicative costs of $n^{O(\log n)}$, or $n^{O(\log n)^{4+\varepsilon}}$, originating in the main subroutines in \S\ref{se:majorou}. However, if we had had at that point a polynomial cost, we could have continued with our bookkeeping until the end and obtained a $n^{O(\log n)}$ runtime: even the $n^{O(\log n)}$ that is weaved already into Cameron's theorem (Theorem~\ref{th:cammar}\eqref{th:cammarsmall}, coming from \cite[Thm.\ 1.1(iii)]{Mar02}) does not pile up eventually, since \eqref{eq:act1} shows that $\max\{e_{1},e_{2}+1\}$ is the correct exponent of the logarithm.

Hence, the bottleneck must be in the subroutines. The local certificates call the algorithm for strings of size $\Omega\left(\frac{n}{\log n}\right)$, for each of the $O(\log n)$-tuples inside an $O(n)$-set: thus, unless one manages to bypass the logarithmic requirement in Lemma~\ref{le:epicfsg}, the routine of Proposition~\ref{pr:computcert} is too expensive to improve the runtime under the $n^{O(\log^{2}n)}$ threshold. Also Split-or-Johnson is in its current form too expensive, but in that case one might make do with reworking the recursion process that comes into play by showing for instance that the worst scenario does not actually happen in real life. It is already a common thread in the literature that distinguishing non-isomorphic graphs is actually pretty easy in general (see \cite{BES80} \cite{BK79}), and a handful of bad cases yields a much worse runtime: SoJ as well analyzes in its recursion hypothetical configurations where it is very difficult to break the symmetry of its vertices, even when we are given from the start that the are few twins among them. It might be feasible to prove that there are actually no such configurations, or alternatively that they are so well-structured that it is possible to describe them entirely and treat them separately as exceptional cases, as was done for instance with the ``three exceptional families'' in \cite[Def.~1.3]{SW16} (the first paper to break the $n^{O(\sqrt{n})}$ threshold on GIP).

\section*{Acknowledgements}

The author thanks H.\ A.\ Helfgott for introducing him to the graph isomorphism problem and for discussions about his paper \cite{Hel19} on the subject.

\bibliography{Bibliography}
\bibliographystyle{alpha}

\begin{mysage}
### Final control strings go here ###
#
# Check whether mlim in {re:thism} is the lowest m for which C(log n)^2<m<=n implies the condition m,n>=mlim for C=Clim.
String_C=''
if RIF(Clim)*loginterval(RIF(mlim))^2>=RIF(mlim) or RIF(Clim)*loginterval(RIF(mlim)-1)^2<RIF(mlim)-1:
    String_C='Error in C in {re:thism}. '
#
# Check whether mlim_free in {re:thism} is the same as Clim_free.
String_Cfree=''
if mlim_free!=Clim_free:
    String_Cfree='Error in C_free in {re:thism}. '
#
# Check whether for eps=1/10 the CFSG-free condition in {re:thism} is stronger than the CFSG one.
epsi=1/10
String_freestrong=''
if RIF(mlim_free)*expinterval(RIF(1/epsi^2))<=RIF(mlim):
    String_freestrong='Error for CFSG-free stronger in {re:thism}. '
#
# Various computational checks in {le:epicfsg}.
String_epicfsg=''
e_np=RIF(5/4)*sqrtinterval(RIF(10000)/loginterval(RIF(10000)))
e_np_up=RIF(e_np*(floor(e_np)+1)/floor(e_np)) #To give a correct integer upper bound.
if 1/2*e_np_up^2*loginterval(e_np_up)*(1+1/loginterval(RIF(2)))>=10000-1:
    String_epicfsg=String_epicfsg+'Error on number of primes in {le:epicfsg}. '
epsi=1/100
if 2/(4+epsi)>=1/2-epsi/10:
    String_epicfsg=String_epicfsg+'Error on epsilon exponent in {le:epicfsg}. '
if expinterval(RIF(1/(10*epsi)))<=RIF(8/5*1/epsi):
    String_epicfsg=String_epicfsg+'Error on cfsg-free bound in {le:epicfsg}. '
#
# Check whether alim_low in {pr:computcert} is compatible with {pr:epicfsgbabai}.
String_alim_low=''
if alim_low<1/loginterval(RIF(2))+2/loginterval(RIF(mlim)):
    String_alim_low='Error in alim_low in {pr:computcert}. '
#
# Check whether alim and alim_low are spaced enough in {pr:computcert}.
String_alim_space=''
if (alim-alim_low)*loginterval(RIF(mlim))<=RIF(1):
    String_alim_space='Error in alim interval in {pr:computcert}. '
# Check whether alim in {pr:computcert} is not higher than m/4.
String_alim=''
if alim>RIF(Clim)*loginterval(RIF(mlim)):
    String_alim='Error in alim in {pr:computcert}. '
#
# Check whether alim_freelow in {pr:computcert} is compatible with {pr:epicfsgbabai}.
String_alim_freelow=''
if alim_freelow<1/loginterval(RIF(2))^5:
    String_alim_freelow='Error in alim_freelow in {pr:computcert}. '
#
# Check whether alim_free and alim_freelow are spaced enough in {pr:computcert}.
epsi=1/10
String_alim_freespace=''
if (alim_free-alim_freelow)*expinterval(1/RIF(epsi))*loginterval(RIF(mlim))^4<=RIF(1):
    String_alim_freespace='Error in alim_free interval in {pr:computcert}. '
# Check whether alim_free in {pr:computcert} is not higher than m/4.
String_alim_free=''
if alim_free>RIF(Clim_free):
    String_alim_free='Error in alim_free in {pr:computcert}. '
#
# Computational check in {pr:computsoj}
String_computsoj=''
if 12/loginterval(RIF(3/2))<=13/loginterval(RIF(2)):
    String_computsoj='Error in {pr:computsoj}. '
#
# Check whether mlim_for_parts in {le:parts} is really the lowest m for which we respect the 2/3 from {Hel19}.
String_m=''
if 1/2*expinterval(2/loginterval(RIF(mlim_for_parts)))>=RIF(2/3) or 1/2*expinterval(2/loginterval(RIF(mlim_for_parts)-1))<RIF(2/3):
    String_m='Error in m in {le:parts}. '
#
# Check whether mlim_for_parts in {le:parts} is really not larger than mlim.
String_m2=''
if mlim_for_parts>mlim:
    String_m2='Error for m too large in {le:parts}. '
#
# Various computational checks in {le:johnjohn}.
String_lejohnjohn=''
if binomial(mplim,2)<11*mplim/2:
    String_lejohnjohn=String_lejohnjohn+'Error for binom(mp,kp)>=11mp/2 in {le:johnjohn}. '
if mplim-4<2/3*mplim:
    String_lejohnjohn=String_lejohnjohn+'Error for 4^(mp-4)>=4^(2/3*mp) in {le:johnjohn}. '
if expinterval(log(3)*4/3*RIF(mplim))/4^4<expinterval(log(4)*2/3*RIF(mplim)):
    String_lejohnjohn=String_lejohnjohn+'Error for 3^(4/3*mp)/4^4>=4^(2/3*mp) in {le:johnjohn}. '
i=mplim
while i<mplim2:
    if binomial(binomial(i,ceil(i/3)),3)<2*binomial(i,floor(i/2)):
        String_lejohnjohn=String_lejohnjohn+'Error for binom(binom(mp,[mp/3]),3)>=2*binom(mp,[mp/2]) in {le:johnjohn}. '
    i=i+1
if expinterval(log(3)*(RIF(mplim2)-3))<=2*expinterval(1/2*log(2*e)*RIF(mplim2)):
    String_lejohnjohn=String_lejohnjohn+'Error for 3^(mp-3)<2*sqrt(2e)^mp in {le:johnjohn}. '
#
# Check whether the exponent in {th:cert} is the correct one.
String_multe_cfsg=''
if 5>=alim_low*loginterval(RIF(mlim)):
    String_multe_cfsg=String_multe_cfsg+'Error for d<=a in {th:cert}. '
if alim+fmsoj>=multe_cfsg:
    String_multe_cfsg=String_multe_cfsg+'Error for multe_cfsg in {th:cert}. '
#
# Various computational checks in {co:cammar}.
String_cocammar=''
if RIF(C0)>=expinterval(log(2)*RIF(Clim)*log(4)^2):
    String_cocammar=String_cocammar+'Error for C0<2^(C*log(4)^2) in {co:cammar}. '
if 1+loginterval(RIF(mlim))/log(2)>=RIF(Clim)*loginterval(RIF(mlim))^2:
    String_cocammar=String_cocammar+'Error for 1+log_2(m)<C*log(n)^2 in {co:cammar}. '
if Clim*log(2)^2<=1:
    String_cocammar=String_cocammar+'Error for k^2*Clog(2)^2>k^2 in {co:cammar}. '
if RIF(Clim)*1/2*loginterval(RIF(mlim)/RIF(Clim))<=loginterval(RIF(mlim)):
    String_cocammar=String_cocammar+'Error for k<sqrt(m/log(m)) in {co:cammar}. '
if RIF(Clim)*loginterval(sqrtinterval(RIF(mlim)/RIF(Clim)))<=loginterval(RIF(mlim)):
    String_cocammar=String_cocammar+'Error for C*log(sqrt(m/C))>log(m) in {co:cammar}. '
if binomial(mplim,floor(mplim/2))>=mlim:
    String_cocammar=String_cocammar+'Error for m=(mp,kp) in {co:cammar}. '
#
# Check whether the 4 sets in {co:certlargesym} are a small expense wrt {th:certsoj}.
String_cocertlargesym=''
if 4>=expinterval(RIF(multe_cfsg)*loginterval(RIF(mlim))^2):
    String_cocertlargesym='Error on 4 sets in {co:certlargesym}. '
#
# Check whether mlim_free1 is the smallest m for which the trivial bound on |H| is not covered by multe_free11*log(m)^3.
String_mlim_free1=''
if (mlim_free1+1/2)*loginterval(RIF(mlim_free1))+1-mlim_free1<=multe_free11*loginterval(RIF(mlim_free1))^3 or (mlim_free1-1+1/2)*loginterval(RIF(mlim_free1)-1)+1-(mlim_free1-1)>multe_free11*loginterval(RIF(mlim_free1)-1)^3:
    String_mlim_free1='Error for the bound on m in {pr:pybersmall}. '
#
# Check whether multe_free12 is the smallest integer above multe_free11 and 51/50*32/log(2)^2.
String_multe_free12=''
if multe_free12<max([multe_free11,51/50*32/log(2)^2]) or multe_free12-1>=max([multe_free11,51/50*32/log(2)^2]):
    String_multe_free12='Error in {pr:pybersmall}. '
#
# Check whether the absorption of bounds in {pr:pybersmall} works out.multe_free1 absorbs multe_free12.
epsi=1/10
String_multe_free1=''
if multe_free12>=Clim_free*expinterval(1/RIF(epsi)^2) or Clim_free!=mlim_free:
    String_multe_free1='Error in bound in {pr:pybersmall}. '
#
# Various computation checks in {pr:pybersmall}.
String_prpybersmall=''
if 4*loginterval(RIF(mlim_free1))/loginterval(RIF(2))<=49.8:
    String_prpybersmall=String_prpybersmall+'Error for 4*log_2(m)>49.8 in {pr:pybersmall}. '
#
# Various computation checks in {pr:pybergiant}.
String_prpybergiant=''
if 3>=Clim_free:
    String_prpybergiant=String_prpybergiant+'Error for d<=a in {pr:pybergiant}. '
if binomial(mplim-1,floor((mplim-1)/2))>=ceil(2/3*mlim):
    String_prpybergiant=String_prpybergiant+'Error for |V_1|<=binom(G,[G/2]) in {pr:pybergiant}. '
epsi=1/10
if alim_free+fmsoj/RIF(1/epsi^2)^3+loginterval(RIF(2))/RIF(1/epsi^2)^5>=multe_free2:
    String_prpybergiant=String_prpybergiant+'Error for multe_free2 in {pr:pybergiant}. '
#
# Check whether the main exponent with CFSG is correct.
String_main_cfsg=''
if main_cfsgtrue>main_cfsg:
    String_main_cfsg=String_main_cfsg+'Error for the main exponent with CFSG. '
if main_cfsgtrue*(1+1/100000)>main_cfsg:
    String_main_cfsg=String_main_cfsg+'Error for the main exponent with CFSG (counting 1+1/100000). '
#
# Check whether the main exponent without CFSG is correct.
String_main_free=''
if main_freetrue>main_free:
    String_main_free=String_main_free+'Error for the main exponent without CFSG. '
if main_freetrue*(1+1/100000)>main_free:
    String_main_free=String_main_free+'Error for the main exponent without CFSG (counting 1+1/100000). '
#
# Various computation checks in the main theorem.
String_main=''
if loginterval(2/3*RIF(mlim))^2>=loginterval(RIF(mlim))^2-3/4*loginterval(RIF(mlim)):
    String_main=String_main+'Error for log(2/3*np)^2<log(np)^2-3/4*log(np). '
if loginterval(RIF(mlim)/(RIF(alim)*loginterval(RIF(mlim))))^2<=1/3*loginterval(RIF(mlim))^2 or loginterval(RIF(mlim))<=1 or 5/3*loginterval(RIF(mlim))^2<=RIF(alim)*loginterval(RIF(mlim)) or 5/3*loginterval(RIF(mlim))^2<=11*loginterval(2*RIF(alim)*loginterval(RIF(mlim))):
    String_main=String_main+'Error for 1st>2nd inside {eq:addt}. '
if loginterval(alim*loginterval(RIF(mlim)))>=1/2*loginterval(RIF(mlim)):
    String_main=String_main+'Error for lognu<1/2*logm. '
if 3*max([alim,alim_free])+2>=2*min([a_cfsg,ab_free/2]) or lowlognu<=2:
    String_main=String_main+'Error for 3*alpha+2<2a. '
if loginterval(2/3*RIF(mlim))^2>=loginterval(RIF(mlim))^2-(3/4+3/100)*loginterval(RIF(mlim)) or expinterval(-3/100*loginterval(RIF(mlim))^2)>=1/4:
    String_main=String_main+'Error for margin of 2nd action. '
if -min([multe_cfsg*loginterval(RIF(mlim))^2,multe_free2*loginterval(RIF(mlim_free))^6])/100000>=loginterval(RIF(49/50)):
    String_main=String_main+'Error for margin of 3rd and 4th action. '
if 49/50+10^(-leewayexp)>=1:
    String_main=String_main+'Error for leeway computation 49/50+10^(-big)<1. '
if ab_free<5 or ab_free!=main_freetrue:
    String_main=String_main+'Error for a,b CFSG-free. '
#
\end{mysage}

$\sage{String_C}\sage{String_Cfree}\sage{String_freestrong}\sage{String_epicfsg}\sage{String_alim_low}\sage{String_alim_space}\sage{String_alim}\sage{String_alim_freelow}\sage{String_alim_freespace}\sage{String_alim_free}\sage{String_computsoj}\sage{String_m}\sage{String_m2}\sage{String_lejohnjohn}\sage{String_multe_cfsg}\sage{String_cocammar}\sage{String_cocertlargesym}\sage{String_mlim_free1}\sage{String_multe_free12}\sage{String_multe_free1}\sage{String_prpybersmall}\sage{String_prpybergiant}
\sage{String_main_cfsg}\sage{String_main_free}\sage{String_main}$

\end{document}